\documentstyle[amscd,amssymb]{amsart}

\newcommand{\seq}[3]{{#1}_{#2}, \ldots, {#1}_{#3}}
\newcommand{\Dws}{D^{(w,\Sigma)}}
\newcommand{\M}{{\cal N}_g}
\newcommand{\Y}{\Sigma \times {{\Bbb S}}^1}
\newcommand{\Spz}{\text{Sp}\, (2g,{\Bbb Z})}
\newcommand{\Spzr}{\text{Sp}\, (2g-2,{\Bbb Z})}
\newcommand{\Cabg}{\CC[\a,\b,\g]}
\newcommand{\hH}{\hat{H}}
\newcommand{\ha}{\hat{\a}}
\newcommand{\hb}{\hat{\b}}
\newcommand{\hg}{\hat{\g}}
\newcommand{\hq}{\hat{\q}}

\newcommand{\he}{\hat{e}}
\newcommand{\Ct}{{\Bbb C}[[t]]}
\newcommand{\Cthabg}{{\Bbb C}[[t]][\hat{\a},\hat{\b},\hat{\g}]}
\newcommand{\whff}{\widetilde{HFF}{}_g^*}
\newcommand{\la}{\langle}
\newcommand{\ra}{\rangle}

\newcommand{\surj}{\twoheadrightarrow}
\newcommand{\inc}{\hookrightarrow}
\newcommand{\ar}{\rightarrow}
\newcommand{\bd}{\partial}
\newcommand{\x}{\times}
\newcommand{\ox}{\otimes}
\newcommand{\iso}{\cong}
\newcommand{\isom}{\stackrel{\simeq}{\ar}}

\newcommand{\point}{\text{pt}}
\newcommand{\CP}{{\Bbb C \Bbb P}}

\newcommand{\Diff}{\text{Diff}}

\newcommand{\rk}{\text{rk}}

\newcommand{\Hom}{\text{Hom}}
\newcommand{\Sym}{\text{Sym}}
\newcommand{\ind}{\text{ind}}
\newcommand{\PD}{\text{P.D.}}

\newcommand{\odd}{\hbox{\scriptsize odd}}
\newcommand{\even}{\hbox{\scriptsize even}}

\newcommand{\id}{\text{id}}

\newcommand{\im}{\text{\rm im\,}}
\newcommand{\red}{\text{red}}

\newcommand{\cF}{{\cal F}}

\newcommand{\cM}{{\cal M}}

\newcommand{\cJ}{{\cal J}}

\newcommand{\cR}{{\cal R}}
\newcommand{\cS}{{\cal S}}

\newcommand{\cV}{{\cal V}}

\renewcommand{\AA}{{\Bbb A}}
\newcommand{\CC}{{\Bbb C}}
\newcommand{\DD}{{\Bbb D}}

\newcommand{\PP}{{\Bbb P}}

\newcommand{\RR}{{\Bbb R}}
\newcommand{\TT}{{\Bbb T}}
\newcommand{\SS}{{\Bbb S}}
\newcommand{\ZZ}{{\Bbb Z}}

\renewcommand{\a}{\alpha}
\renewcommand{\b}{\beta}
\renewcommand{\d}{\delta}
\newcommand{\g}{\gamma}
\newcommand{\e}{\varepsilon}

\newcommand{\s}{\sigma}

\newcommand{\p}{\phi}
\newcommand{\q}{\psi}
\renewcommand{\t}{\tau}
\newcommand{\r}{\rho}
\renewcommand{\v}{\varphi}

\renewcommand{\S}{\Sigma}
\newcommand{\D}{\Delta}
\renewcommand{\L}{\Lambda}
\renewcommand{\P}{\Phi}

\theoremstyle{plain}
\begingroup
\theorembodyfont{\sl}
\newtheorem{thm}{Theorem}[section]
\newtheorem{cor}[thm]{Corollary}
\newtheorem{lem}[thm]{Lemma}
\newtheorem{prop}[thm]{Proposition}

\endgroup

\theoremstyle{definition}
\newtheorem{defn}[thm]{Definition}

\theoremstyle{remark}
\newtheorem{rem}[thm]{Remark}

\numberwithin{equation}{section}

\title{Fukaya-Floer homology of $\Y$ and Applications}

\author{Vicente Mu\~noz}
\address{Departamento de \'Algebra, Geometr\'{\i}a y Topolog\'{\i}a \\ Facultad de Ciencias \\ 
Universidad de M\'alaga \\ 29071 M\'alaga \\ Spain}
\email{vmunoz@@agt.cie.uma.es}
\thanks{\hbox{$^*$}Supported by a grant from Ministerio de Educaci\'on y Cultura of Spain \\
Key words: Fukaya-Floer homology, Floer homology, $4$-manifolds,
Donaldson invariantes, simple type. \\
Mathematical Subject Classification. Primary: 58D27. Secondary: 57R57.}
\date{March, 1998. Revised April, 1999.}

\begin{document}

\maketitle

\begin{abstract}
  We determine the Fukaya-Floer (co)homology groups of the three-mani\-fold $Y=\Y$, where
  $\S$ is a Riemann surface of genus $g \geq 1$. These are of two kinds. For the $1$-cycle
  $\SS^1 \subset Y$, we compute the Fukaya-Floer cohomology $HFF^*(Y,\SS^1)$ and
  its ring structure, which is a sort
  of deformation of the Floer cohomology $HF^*(Y)$. On the other hand, for $1$-cycles
  $\d \subset \S \subset Y$, we determine the Fukaya-Floer homology
  $HFF_*(Y,\d)$ and its $HF^*(Y)$-module structure.

  \noindent We give the following applications: \\
  \noindent $\bullet$~We show that every four-manifold with $b^+>1$ is of finite type. \\
  \noindent $\bullet$~Four-manifolds which arise as connected sums along surfaces of four-manifolds with 
  $b_1=0$ are of simple type and we give constraints on their basic classes. \\
  \noindent $\bullet$~We find the invariants of the product of two Riemann surfaces both 
  of genus greater or equal than one.
\end{abstract}

\section{Introduction}
\label{sec:f1}

The structure of Donaldson invariants of $4$-manifolds has been found out by 
Kronheimer and Mrowka~\cite{KM} and Fintushel and Stern~\cite{FS} for a large
class of $4$-manifolds (those of simple type with $b_1=0$, $b^+>1$) making use of universal
relations coming from embedded surfaces. In order to analyse general $4$-manifolds, we
need to set up first the right framework for getting enough universal relations. 
It is the purpose of this work to do this by using the Fukaya-Floer homology of the
three manifold $Y=\Y$, the product of a surface times a circle. This is obviously not
the only way, but it already gives new results.

Donaldson invariants for a $4$-manifold $X$ with $b^+>1$ are defined as linear functionals
$$
  D^w_X: \AA(X)= \Sym^*(H_0(X) \oplus H_2(X)) \ox \L^* H_1(X) \ar \CC,
$$
where $w \in H^2(X;\ZZ)$. For the homology $H_*(X)$ we shall understand complex coefficients.
$\AA(X)$ is graded giving degree $4-i$ to the elements in $H_i(X)$. 
There is a slight difference in our definition of $\AA(X)$ with that of
Kronheimer and Mrowka~\cite{KM},
as we do not consider $3$-homology classes (this is done in this way
since the techniques here contained are
not well suited to deal with these classes).

We say that $X$ is of {\bf $w$-simple type} when $D^w_X((x^2-4)z)=0$, for any $z\in \AA(X)$.
If $X$ has $b_1=0$ and it 
is of $w$-simple type, then it is of $w'$-simple type, for any other $w'$ and it is
said to be of simple type for brevity.
Analogously, we say that $X$ is of {\bf $w$-finite type} when there 
is some $n \geq 0$ such that 
$D^w_X((x^2-4)^nz)=0$, for any $z\in \AA(X)$. The order is the minimum such $n$, so
order $1$ means simple type and order $0$ means that the Donaldson invariants are identically zero.
$X$ is of finite type if it is of $w$-finite type for any $w$.
For completeness we introduce the notion of $X$ being of {\bf $w$-strong simple type} when 
$D^w_X((x^2-4)z)=0$, for any $z\in \AA(X)$ and $D^w_X(\g z)=0$, for any $\g \in H_1(X)$ and
any $z\in \AA(X)$. This condition is the right one for extending the concept of
simple type in the case
$b_1=0$ to the case $b_1>0$. It gives the same structure theorem for the invariants
presented in~\cite{KM}.
Introducing $3$-homology classes, we would have (potentially) different definitions,
but we shall not deal here with that issue.
Also the order of $w$-finite type of $X$ does not depend on $w$ (see~\cite{otro}).

Let $\S=\S_g$ be a Riemann surface of genus $g\geq 1$ and consider the three-manifold
$Y=\Y$. In~\cite{floer} we computed the ring structure of the (instanton) Floer (co)homology of
$Y$, together with the $SO(3)$-bundle with $w_2=\PD[\SS^1]$.
This gadget encondes all the relations
$R\in \AA(\S)$ satisfied by all $4$-manifolds $X$ containing an embedded surface $\S$, 
representing an odd homology class and with $\S^2=0$. More accurately, for such $X$, 
$D^w_X(R z)=0$, for any $z\in \AA(\S^{\perp})$ and $w \in H^2(X;\ZZ)$ with $w \cdot\S
\equiv 1\pmod 2$. This 
is so since we have a decomposition $X=X_1\cup_Y A$, where $A$ is a
tubular neighbourhood of $\S$, and we can consider $R \in \AA(A)$ and $z \in \AA(X_1)$. Then
the (relative) Donaldson invariants for $A$ corresponding to $R$ are already vanishing.

In order to drop the condition $z\in \AA(\S^{\perp})$, the useful space to 
work in is no longer the Floer homology, but the 
extension developed by Fukaya~\cite{Fukaya}~\cite{HFF} and known as
Fukaya-Floer homology. This one deals with $2$-cycles in $X$ cutting $Y$ non-trivially. In our case as
$X=X_1\cup_Y A$, the only possibility for the cutting of a $2$-cycle of $X$
with $Y$ is $n\SS^1 \subset Y=\Y$.
Here we extend the arguments of~\cite{floer}
to find (to a large extent) the structure of the
Fukaya-Floer (co)homology $HFF^*(Y,n\SS^1)$ (with the same $SO(3)$-bundle as above). 

This will set up the background work
necessary to give a structure theorem of the Donaldson invariants
for manifolds not of simple type~\cite{Kr}, work which will be carried out 
in future. Such a structure theorem was conjectured in~\cite{Kr} and 
presumably, it might follow from the arguments given in~\cite{KM}~\cite{FS}.
The first result in this direction is the finite type condition for all 
$4$-manifolds with $b^+>1$, which we prove. 
Fr{\o}yshov~\cite{froyshov} and Wieczorek~\cite{wieczorek} 
have given alternative proofs only valid for simply connected $4$-manifolds.

On the other hand, there is another possibility for the Fukaya-Floer homology of $Y$, 
that is cutting with $\d \subset \S \subset Y$, $\d$ primitive in homology. For completeness, 
we also determine the structure of 
$HFF_*(Y,\d)$ and, as an application we show that a connected sum
of two $4$-manifolds with $b_1=0$ along a surface (which represents an odd homology class and has 
self-intersection zero) is of simple type and give serious constraints on its basic classes,
along the lines in~\cite{genus2}~\cite{genusg}.

Finally, we prove that the product of two Riemann surfaces are of strong simple type using 
these techniques and give their Donaldson invariants. The basic classes coincide with its Seiberg-Witten
basic classes, as expected.

The paper is organised as follows. In sections~\ref{sec:f2} and~\ref{sec:f3}
we review the construction of the Floer homology and Fukaya-Floer homology 
of a three manifold with $b_1 \neq 0$. Then in section~\ref{sec:f4} we recall, 
for the convenience of the reader, the structure of the Floer cohomology 
$HF^*(\Y)$.
Section~\ref{sec:f5} is devoted to determining the Fukaya-Floer cohomology 
corresponding to the $1$-cycle $\SS^1\subset Y$, 
$HFF^*(\Y,\SS^1)$, which is finally given in terms of the relations 
(which are partially determined) satisfied by the generators.
On the other hand, the Fukaya-Floer cohomology $HFF^*(\Y,\d)$ corresponding 
to $\d \subset \S \subset \Y$, is given in 
section~\ref{sec:f6}. The applications already mentioned are collected in section~\ref{sec:f7}.

\section{Review of Floer homology}
\label{sec:f2}

In this section we are going to review the construction of the Floer homology 
groups of a $3$-manifold $Y$ with $b_1 >0$, endowed with an $SO(3)$-bundle
$P$ with second Stiefel-Whitney class $w_2
=w_2(P) \neq 0 \in H^2(Y;\ZZ/2\ZZ)$. Recall that $w_2$ determines uniquely $P$.
To be more precise, we are going
to suppose that $w_2$ has an integral lift (i.e.\ that the $SO(3)$-bundle lifts
to an $U(2)$-bundle). This case is in contrast with the 
case of Floer homology of rational homology spheres. All the facts
stated here are well known. For full treatment and
proofs see~\cite{Fl}~\cite{Don1}~\cite{Furuta} (for the case
of Floer homology of rational homology spheres see~\cite{Braam}).
We shall use complex coefficients for the Floer homology (although it
is usually developed over the integers).

\subsection{Floer homology}
As $w_2 \neq 0 \in H^2(Y;\ZZ/2\ZZ)$, 
there are no reducible flat connections on $P$. We say
that $P$ is free of flat reductions. Possibly after a small
perturbation of the flat equations, there will be finitely many flat
connections $\r_j$, and they will all be non-degenerate. The Floer
complex $CF_*(Y)$ is the complex vector space with basis given by the $\r_j$,
with a $\ZZ/4\ZZ$-grading which is given by the index~\cite{Don1}~\cite{Don2}.
Actually, this grading is only defined up to addition of a constant.
The complex $CF_*(Y)$ depends on
$w_2$, but in general we will not express this in the notation. 

Now for every two flat connections $\r_k$ and $\r_l$ there is a moduli
space $\cM(\r_k, \r_l)$ of (perturbed) ASD connections on the tube
$Y \x \RR$ with limits $\r_k$ and $\r_l$. There is an $\RR$-action by traslations
and $\cM_0(\r_k, \r_l)$ shall stand for the quotient 
$\cM(\r_k, \r_l)/\RR$. This space has
components $\cM_0^D(\r_k, \r_l)$ of dimensions $D \equiv 
\ind(\r_k)-\ind(\r_l)-1 \pmod 4$, 
and can be oriented in a compatible way~\cite{Fl}. When 
$\ind(\r_l)=\ind(\r_k)-1$, there is a compact zero dimensional moduli
space $\cM_0^0(\r_k, \r_l)$, for which the algebraic 
number of its points is defined $\# \cM_0^0(\r_k, \r_l)$. The
boundary map of the Floer complex is then
\begin{eqnarray*}
  \bd:  CF_i(Y) & \ar & CF_{i-1}(Y) \\ 
        \r_k & \mapsto &  \hspace{-5mm}
  \sum_{\r_l \atop \ind(\r_l)=\ind(\r_k)-1}\hspace{-5mm} 
  \# \cM_0^0(\r_k, \r_l) \r_l
\end{eqnarray*}

To see that $(CF_*(Y), \bd)$ defines a complex we need to check that $\bd^2=0$ 
(see~\cite{Don1}~\cite{Furuta}). For that, 
consider $\r_k$ and $\r_l$ flat connections such that
  $\ind(\r_l)=\ind(\r_k)-2$. Then the moduli space $\cM_0^1(\r_k, \r_l)$
  is a smooth one dimensional manifold which can be compactified adding
  the broken instantons in 
  \begin{equation}
    \bigcup_{\r_m \atop
    \ind(\r_m)=\ind(\r_k)-1}\hspace{-5mm} \cM_0^0(\r_k, \r_m) \x
    \cM_0^0(\r_m, \r_l). \label{eqn:f2.1}  
  \end{equation}
So this compactification, $\overline{\cM}_0^1(\r_k, \r_l)$, is a
manifold with boundary given by~\eqref{eqn:f2.1}. The homology class of the
boundary is zero, i.e.
  $$
    \sum_{\r_m \atop
    \ind(\r_m)=\ind(\r_k)-1}\hspace{-5mm} \# \cM_0^0(\r_k, \r_m) \cdot
    \# \cM_0^0(\r_m, \r_l) =0,
  $$
or equivalently, $\bd^2=0$.  

We define the
Floer homology $HF_*(Y)$ as the homology of this complex $(CF_*(Y),\bd)$ 
(see~\cite{Fl}).
It can be proved that these groups do not
depend on the metric of $Y$ or on the chosen perturbation of the ASD
equations. The groups $HF_*(Y)$ are natural under diffeomorphisms of
the pair $(Y,P)$.
The Floer cohomology $HF^*(Y)$ is defined analogously out of
the dual complex $CF^*(Y)$ and it is naturally isomorphic to $HF_{c-*}(\overline
Y)$, where $\overline Y$ denotes $Y$ with reversed
orientation ($c$ is a constant that we need to introduce due to the indeterminacy 
of the grading). 
The natural pairing $HF_*(Y) \ox HF^*(Y) \to \CC$ yields
the pairing
$\la,\ra: HF_*(Y) \ox HF_{c-*} (\overline Y) \to \CC$.
It is worth noticing that when 
$Y$ has an orientation reversing diffeomorphism, i.e. $Y\iso
\overline Y$, we have a pairing
\begin{equation}
  \la,\ra: HF_*(Y) \ox HF_{c-*} (Y) \to \CC.
  \label{eqn:f2.2}
\end{equation}

\subsection{Action of $H_*(Y)$ on $HF_*(Y)$}
Let $\a \in H_{3-i}(Y)$. We have cycles $V_{\a}$, in the moduli spaces
$\cM(\r_k, \r_l)$, of codimension $i+1$, representing $\mu(\a \x
\point)$, for $\a\x\point \subset Y\x\RR$,
much in the same way as in the case of
a closed manifold~\cite{DK}~\cite{KM}. Using them, we construct a map
\begin{eqnarray*}
  \mu(\a):  CF_j(Y) & \ar & CF_{j-i-1}(Y) \\ 
          \r_k & \mapsto & \hspace{-5mm}
  \sum_{\r_l \atop \ind(\r_l)=\ind(\r_k)-i-1} \hspace{-5mm}(\#
  \cM^{i+1}(\r_k, \r_l) 
      \cap V_{\a}) \, \r_l 
\end{eqnarray*}
(note that this time we do not quotient by the traslations
as the cycles $V_{\a}$ are not translation invariant).
For $\ind(\r_l)=\ind(\r_k)-i-2$ consider the $1$-dimensional moduli 
space $\cM^{i+2}(\r_k, \r_l) \cap V_{\a}$, and 
looking at the number of points in the boundary of its compactification,
as we did before, we get that
$$
  \sum_{\r_s \atop \ind(\r_s)=\ind(\r_l)+1} \hspace{-7mm} (\#
  \cM^{i+1}(\r_k, \r_s)  \cap
  V_{\a}) \cdot \# \cM^0_0(\r_s,\r_l) 
  + \hspace{-5mm} \sum_{\r_s \atop \ind(\r_s)=\ind(\r_k)-1} \hspace{-7mm} \#
  \cM^0_0(\r_k, \r_s) \cdot ( \# \cM^{i+1}(\r_s,\r_l) \cap V_{\a}) 
$$
vanishes, or equivalently, $\bd \circ \mu(\a) +\mu(\a)\circ \bd=0$.
So $\mu(\a)$ descends to a map
$$
    \mu(\a) : HF_*(Y) \to HF_{*-i-1}(Y).
$$

\subsection{Products in Floer homology}
Another useful construction, which is outlined in~\cite{Don1},
is the following.
Suppose that we have an (oriented) four-dimensional cobordism $X$ between 
two closed oriented $3$-manifolds $Y_1$ and $Y_2$, i.e. $X$ is a $4$-manifold with
boundary $\bd X=Y_1 \sqcup \overline Y_2$. Suppose that we have an
$SO(3)$-bundle $P_X$ over $X$ such that $P_1=P_X|_{Y_1}$ and $P_2=P_X|_{Y_2}$
satisfy $w_2(P_i)\neq 0$, $i=1,2$, so that we have defined the 
Floer homologies of $(Y_1,P_1)$ and $(Y_2,P_2)$. 
Furnishing $X$ with two cylindrical ends,
the cobordism $X$ gives a map
\begin{eqnarray*}
  \P_X: CF_*(Y_1) &\ar & CF_*(Y_2) \\
           \r_k & \mapsto & \sum_{\r'_l} \# \cM^0(X,\r_k, \r'_l) \r'_l,
\end{eqnarray*}
where $\cM(X,\r_k, \r'_l)$ is the moduli space of (perturbed) ASD connections
on $X$ with flat limits $\r_k$ on the $Y_1$ side and $\r'_l$ on the $Y_2$ side.
We can check along the lines above that $\bd \circ \P_X + \P_X \circ \bd =0$,
so that $\P_X$ defines a map $\P_X:HF_*(Y_1) \ar HF_*(Y_2)$. 
Also note that if $\a_1 \in H_*(Y_1)$ and $\a_2 \in H_*(Y_2)$ define the same homology class
in $X$, then $\mu(\a_2) \circ \P_X=\P_X\circ \mu(\a_1)$. 

On the other hand, suppose that $Y_1$ and $Y_2$ are oriented $3$-manifolds and
$P_1$ and $P_2$ are $SO(3)$-bundles with $w_2(P_i)\neq 0$, $i=1,2$.
Consider the $SO(3)$-bundle $P=P_1\sqcup P_2$ over $Y=Y_1 \sqcup Y_2$. Every flat
connection on $P$ is of the form $\t_{ll'}=(\r_l,\r'_{l'})$ and $\ind (\t_{ll'})=
\ind (\r_l) +\ind(\r'_{l'})$. So we have naturally
$CF_*(Y)=CF_*(Y_1) \ox CF_*(Y_2)$. For $\ind (\t_{kk'})=\ind (\t_{ll'})+1$,
if $\cM^0_0(\t_{kk'},\t_{ll'})$ is not empty, then either $\r_k=\r_l$ and
$$
   \cM^0_0(\t_{kk'},\t_{ll'})=\cM^0_0((\r_k,\r'_{k'}),(\r_l,\r'_{l'}))= \{\r_k\} \x 
   \cM^0_0(\r'_{k'},\r'_{l'})
$$
or $\r'_{l'}=\r'_{k'}$ and $\cM^0_0(\t_{kk'},\t_{ll'})=\cM^0_0(\r_k,\r_l) \x\{\r'_{k'}\}$.
This proves that 
$\bd_{CF_*(Y)}=\bd_{CF_*(Y_1)}+\bd_{CF_*(Y_2)}$ and therefore 
$HF_*(Y)= HF_*(Y_1) \ox HF_*(Y_2)$.

Putting the above together, a product for $HF_*(Y)$ might arise 
as follows. Suppose that there is a cobordism
between $Y \sqcup \, Y$ and $Y$, i.e. an oriented $4$-manifold $X$ with
boundary $\bd X= Y\sqcup Y \sqcup \overline Y$. 
Then there is a map 
$$
  HF_*(Y) \ox HF_*(Y) \ar HF_*(Y).
$$
In some particular cases, 
this gives an associative and graded commutative ring structure on $HF_*(Y)$.
We shall prove it for the particular $3$-manifold $Y=\Y$ using an argument along 
different lines (see section~\ref{sec:f4}).

\subsection{Relative invariants of $4$-manifolds}
Our purpose now is to recall the definition of Donaldson invariants of 
an (oriented) $4$-manifold $X$ with boundary $\bd
X=Y$, for any 
$w \in H^2(X;\ZZ)$ such that $w|_Y=w_2 \in H^2(Y;\ZZ/2\ZZ)$. These invariants will not
be numerical (in contrast with the case of a closed $4$-manifold),
instead they live in the Floer homology $HF_*(Y)$.

We give $X$ a cylindrical end modelled on $Y \x [0,\infty)$ and denote
it by $X$ again (no confusion should arise out of this).
We have moduli spaces $\cM(X, \r_l)$ of (perturbed) ASD connections with
finite action and asymptotic to $\r_l$. $\cM(X, \r_l)$ has components
$\cM^D(X, \r_l)$ of dimensions $D \equiv
\ind(\r_l)+ C \pmod 4$, for some fixed constant $C$ only dependent on $X$.
The spaces $\cM(X, \r_l)$ can be oriented coherently and,
for $z=\a_1 \a_2 \cdots \a_r \in \AA(X)$ of degree $d$,
we can choose (generic) cycles $V_{\a_i}
\subset \cM(X, \r_l)$ representing $\mu(\a_i)$, so that we have defined an element
$$
  \p^w(X,z)=\sum_{\r_l \atop \ind(\r_l)+C=d} (\# \cM^d(X, \r_l)
      \cap V_{\a_1} \cap \cdots \cap V_{\a_r}) \,\r_l  \in CF_*(Y).
$$
This element has boundary zero (we recommend the reader prove this fact)
and hence it defines a homology class, which is called the {\bf relative
invariants} of $X$, denoted again by $\p^w(X,z)$, in
$HF_*(Y)$ (see~\cite{Don1}~\cite{Furuta}). 

Once defined the relative invariants, we have a gluing theorem for them. 
Suppose we are in the situation of a closed $4$-manifold 
$X=X_1 \cup_Y X_2$, obtained as the union of two $4$-manifolds with boundary,
where $\bd X_1=Y$ and $\bd X_2=\overline Y$. Let 
$w \in H^2(X; \ZZ)$ with $w|_Y=w_2 \in H^2(Y;\ZZ/2\ZZ)$ as above (this implies 
in particular $b^+(X)>0$, so the Donaldson invariants of $X$ are defined; in the case
$b^+=1$ relative to chambers~\cite{Kotschick}~\cite{wall}). 
We need another bit of terminology from~\cite{genus2}

\begin{defn}
\label{def:f2.allowable}
  $(w,\S)$ is an {\bf allowable} pair if 
  $w, \S \in H^2(X; \ZZ)$, $w \cdot \S \equiv 1\pmod 2$ and $\S^2 =0$. 
  Then we define $D^{(w,\S)}_X=D^w_X +D^{w+\S}_X$.
\end{defn}

  Usually, for $X=X_1\cup_Y X_2$, we have $w \in
  H^2(X;\ZZ)$ with $w|_Y=w_2$ as above, and 
  $\S \in H^2(X;\ZZ)$ whose Poincar\'e dual lies in the image of
  $H_2(Y;\ZZ) \ar H_2(X;\ZZ)$,
  and satisfies $w \cdot \S \equiv 1\pmod 2$. Then $(w,\S)$ is an allowable pair.
The series $\Dws_X$ behaves much in the same way as the Kronheimer-Mrowka~\cite{KM} series
$\DD^w_X(\a)=D^w_X((1+{x\over 2})e^{\a})$ (they are equivalent for 
manifolds of simple type with $b_1=0$ and $b^+>1$, see~\cite{genusg} for an
explicit formula), but it is a more efficient way of wrapping around the
information in general. 

When $b^+=1$, the Donaldson invariants depend on
the choice of metric for $X$. In general, we shall consider a family of metrics $g_R$, $R>1$,
giving a neck of length $R$, i.e. $X=X_1 \cup (Y\x [0,R]) \cup X_2$, where the metrics
on $X_1$ and $X_2$ are fixed and the metric on $Y\x [0,R]$ is of the form $g_Y +dt\ox dt$, for
a fixed metric $g_Y$ on $Y$.
Then for large enough $R$ (depending on the degree of $z\in\AA(X)$), the metrics $g_R$ stay 
within a fixed chamber and $D^w_X(z)$ is well deined. We shall refer to these metrics as
{\em metrics on $X$ giving a long neck}. Note that in this case $\p^w(X_i,z_i)$ also depends
on the metric on $X_i$.

\begin{thm}
\label{thm:f2.Floer}
  Let $X=X_1 \cup_Y X_2$ be as above and $w \in
  H^2(X;\ZZ)$ with $w|_Y=w_2$. Take 
  $\S \in H^2(X;\ZZ)$ whose Poincar\'e dual lies in the image of
  $H_2(Y;\ZZ) \ar H_2(X;\ZZ)$,
  and satisfies $w \cdot \S \equiv 1\pmod 2$. Put
  $w_i =w|_{X_i} \in H^2(X_i;\ZZ)$.
  For $z_i \in \AA(X_i)$, $i=1,2$, it is
  $$
     \Dws_X(z_1\,z_2)=\la\p^{w_1}(X_1,z_1),
     \p^{w_2}(X_2,z_2)\ra.
  $$
  When $b^+=1$ the invariants are calculated for metrics on $X$ giving a 
  long neck.
\end{thm}

This is a standard and well known fact~\cite{Don2}. The only not-so-standard fact
is the appearance of $(w,\S)$. This is so since we are working
with $SO(3)$-Floer theory instead of $U(2)$-Floer theory which would give Floer groups 
graded modulo $8$. When we glue the $SO(3)$-bundles over $X_1$ and $X_2$
with second Stiefel-Whitney classes $w_1$ and $w_2$ we can do it in different 
ways, as there is a choice of gluing automorphism of the bundles along $Y$, 
and both $w$ and $w+\S$ are two different possibilities for the resulting
$SO(3)$-bundle whose difference in the indices of both is $4$
(see~\cite{Don1}~\cite{Don2}).

In general we shall write
$$ 
  \p^w(X,e^{t\a}) = \sum_d {\p^w(X,\a^d) \over d!}t^d,
$$
as an element living in $HF_*(Y)\ox \Ct$. Theorem~\ref{thm:f2.Floer}
can be rewritten as

\begin{thm}
\label{thm:f2.Floer.series}
  Let $X=X_1 \cup_Y X_2$ be as above and $w \in
  H^2(X;\ZZ)$ with $w|_Y=w_2$. Take 
  $\S \in H^2(X;\ZZ)$ whose Poincar\'e dual lies in the image of
  $H_2(Y;\ZZ) \ar H_2(X;\ZZ)$,
  and satisfies $w \cdot \S \equiv 1\pmod 2$. Put
  $w_i =w|_{X_i} \in H^2(X_i;\ZZ)$. Then for $\a_i \in H_2(X_i)$, $i=1,2$,
  $$
     \Dws_X(e^{t(\a_1+\a_2)}) 
     =\la\p^{w_1}(X_1,e^{t\a_1}),\p^{w_2}(X_2,e^{t\a_2})\ra.
  $$
  When $b^+=1$ the invariants are calculated for metrics on $X$ giving a 
  long neck.
\end{thm}

\section{Review of Fukaya-Floer homology}
\label{sec:f3}

Now we pass on to the definition of
the Fukaya-Floer homology groups, which are a refinement of the Floer homology 
groups of a $3$-manifold $Y$ with $b_1>0$. The construction is
initially given by Fukaya in~\cite{HFF}
and explained in a paper worth reading by Braam and 
Donaldson~\cite{Fukaya}.
The origin of the Fukaya-Floer homology is the need of defining relative 
invariants (and establishing the appropriate gluing threorem) for $2$-homology classes
crossing the neck in a splitting $X=X_1\cup_Y X_2$. They are in some sense more
natural than the Floer homology
from the point of view of the Donaldson invariants of $4$-manifolds. 

\subsection{Fukaya-Floer homology}
The input is a triple $(Y, P, \d)$,
where $P$ is an $SO(3)$-bundle with $w_2 \neq 0$ over an oriented $3$-manifold $Y$ and
$\d$ is a loop in $Y$, i.e.\ an (oriented) embedded $\d \iso \SS^1 \inc Y$.
The complex $CFF_*(Y,\d)$ will be the total complex of the double complex
$$
CF_*(Y) \otimes \hH_*(\CP^{\infty}),
$$
where $\hH_*(\CP^{\infty})$ is the
completion of $H_*(\CP^{\infty})$, i.e.\ the ring of formal power series. 
Recall that $H_i(\CP^{\infty}) = 0$
for $i$ odd and $\CC$ for $i$ even (we are using complex coefficients). Therefore
\begin{equation}
  CFF_i(Y,\d)=CF_i(Y)\x CF_{i-2}(Y) t \x CF_{i-4}(Y) {t^2 \over 2!} \x 
  CF_{i-6}(Y) {t^3 \over 3!} \x \cdots
\label{eqn:f3.1}
\end{equation}
The labels ${t^k \over k!}$ must be understood as the generators of $H_{2k}
(\CP^{\infty})$ and have an assigned (homological) degree $2k$.
So $CFF_*(Y,\d)=CF_*(Y) \ox \Ct$, i.e.
Fukaya-Floer chains are infinite sequences of (possibly non-zero) Floer chains.
This complex is also graded over $\ZZ/4\ZZ$. 
To construct the boundary $\bd$ we work as follows.
For every pair of flat connections $\r_k$ and $\r_l$ 
we have the moduli space $\cM_0(\r_k,\r_l)$ of section~\ref{sec:f2} and we can
construct generic cycles $V_{\d \x\RR}$ representing $\mu(\d\x\RR)$, 
for $\d \x\RR\subset Y\x\RR$, which intersect transversely in the 
top stratum of the compactification of $\cM_0(\r_k, \r_l)$.
The boundary of $CFF_*(Y)$ will be defined as (see~\cite{Fukaya})
\begin{eqnarray*}
  \bd:  CFF_i(Y) & \ar & CFF_{i-1}(Y) \\
        \r_k {t^a \over a!}& \mapsto & \sum_{\r_l\atop b\geq a} {b \choose a}
       (\# \cM^{2(b-a)}_0(\r_k, \r_l)\cap V_{\d \x\RR}^{b-a}) \r_l {t^b \over b!}
\end{eqnarray*}
for $\r_k \in CF_{i-2a}$, $\r_l \in CF_{i-1-2b}$. Here $V_{\d \x\RR}^{b-a}$
means the intersection of $b-a$ different generic cycles (we only have added
the labels to the formula in~\cite{Fukaya}).
The proof of $\bd^2=0$ is given in~\cite{Fukaya} and runs as follows.
  Consider two flat connections $\r_k$ and $\r_l$,
  such that $\ind(\r_l)=\ind(\r_k)-2-2e$. Then 
  the moduli space $\cM^{2e+1}_0(\r_k, \r_l) \cap V_{\d\x\RR}^e$
  is a one dimensional manifold. We compactify it and count the
  boundary points in the same way as in the case of Floer homology to get
  $$
    \sum_{\r_m \atop
    \ind(\r_m)=\ind(\r_k)-1-2f}\hspace{-5mm} {e \choose f}
    \# \cM^{2f}_0(\r_k, \r_m) \cap V_{\d\x\RR}^f\cdot
    \# \cM^{2(e-f)}_0(\r_m, \r_l) \cap V_{\d\x\RR}^{e-f} =0,
  $$
  equivalently $\bd^2 \r_k= 0$.

We have thus defined the Fukaya-Floer homology $HFF_*(Y,\d)$
as the homology of the complex $(CFF_*(Y,\d),\bd)$. These groups are
independent of metrics and of perturbations of equations~\cite{HFF}.
For the effective computation of $HFF_*(Y,\d)$, we construct a spectral
sequence next.
There is a filtration $(K^{(i)})_* = CF_*(Y) \otimes (\prod_{* \geq i}
\hH_*(\CP^{\infty}))$ of $CFF_*(Y,\d)$ inducing a spectral sequence whose
$E_3$ term is $HF_*(Y) \otimes \hH_*(\CP^{\infty})$ and converging to
the Fukaya-Floer groups (there is no problem of convergence because of
the periodicity of the spectral sequence). The boundary $d_3$
turns out to be
$$
   \mu(\d): HF_i(Y) \otimes H_{2j}(\CP^{\infty}) \ar HF_{i-3}(Y)
   \otimes H_{2j+2} (\CP^{\infty}).
$$

The obvious $\Ct$-module structure of $CFF_*(Y,\d)=
CF_*(Y)\ox\Ct$ descends to give a 
$\Ct$-module structure for $HFF_*(Y,\d)$ (the boundary $\bd$ is $\Ct$-linear
thanks to the choice of denominators
in~\eqref{eqn:f3.1}).

The Fukaya-Floer cohomology will be defined as the homology of the
dual complex $CFF^*(Y,\d)=\Hom_{\Ct}(CFF_*(Y,\d),\Ct)$. We remark that this
is a different definition from that of~\cite{Fukaya}. There is a pairing
$\la,\ra:HFF_*(Y,\d) \otimes HFF^*(Y, \d) \ar \Ct$ and an isomorphism
$HFF_*(\overline Y, -\d) \iso HFF^*(Y,\d)$,
where $-\d$ is $\d$ with reversed orientation, hence 
a pairing for the Fukaya-Floer homology groups 
$$ 
   \la,\ra:HFF_*(Y,\d) \otimes HFF_*(\overline Y, -\d) \ar \Ct.
$$
This can be defined through the spectral sequence from
the natural pairing in $HF_*(Y)$. Also it is a nice way
of collecting all the pairings $\s_m$ in~\cite{Fukaya}. 

The Fukaya-Floer homology may also be defined for $(Y,P,\d)$ where $\d \iso \SS^1 \sqcup
\cdots \sqcup \SS^1 \inc Y$ is a collection of finitely many disjoint loops (possibly none).
In particular, for $\d=\o$, $HFF_*(Y,\o)=HF_*(Y)\ox \Ct$ naturally.

\subsection{Action of $H_*(Y)$ on $HFF_*(Y,\d)$}
This is explained in~\cite[section 5.3]{thesis}.
Let $\a \in H_{3-i}(Y)$. We define $\mu(\a)$ at the level of chains as
\begin{eqnarray*}
  \mu(\a):  CFF_j(Y) & \ar & CFF_{j-i-1}(Y) \\ 
          \r_k {t^a \over a!} & \mapsto & \sum_{\r_l\atop b\geq a} {b \choose a}
       (\# \cM^{2(b-a)+i+1}(\r_k, \r_l)\cap V_{\d \x\RR}^{b-a} \cap V_{\a \x \point}) 
      \r_l {t^b \over b!}
\end{eqnarray*}
for $\r_k \in CF_{j-2a}$, $\r_l \in CF_{j-i-1-2b}$. 
Again $\bd \circ \mu(\a) +\mu(\a)\circ \bd=0$
and $\mu(\a)$ descends to a map
$\mu(\a):HFF_*(Y,\d)\ar HFF_{*-i-1}(Y,\d)$. 
For instance, for $HFF_*(Y,\o)=HF_*(Y)\ox \Ct$, the map $\mu(\a)$ is the one induced 
from $HF_*(Y)$. In general, the induced map
in the term $E_3=HF_*(Y)\ox \Ct$ of the spectral sequence computing 
$HFF_*(Y,\d)$ is $\mu(\a)$ in Floer homology.
The structure of the map $\mu(\a)$ is the cornerstone of the analysis 
in~\cite[chapter 5]{thesis} and the seed of this work. 

\subsection{Products in Fukaya-Floer homology}
We can extend the arguments of section~\ref{sec:f2}.
Suppose that we have an (oriented) four-dimensional cobordism $(X,D,P)$ between 
two triples $(Y_1,\d_1,P_1)$ and $(Y_2,\d_2,P_2)$ as above. Then $\P_X$ is defined at
the level of chains by
\begin{eqnarray*}
  \P_X: CFF_*(Y_1) &\ar & CFF_*(Y_2) \\
           \r_k {t^a \over a!} & \mapsto & 
   \sum_{\r'_l} {b \choose a} (\# \cM^{2(b-a)}(X,\r_k, \r'_l) \cap V_D^{b-a} )
   \r'_l {t^b \over b!}.
\end{eqnarray*}
As $\bd \circ \P_X + \P_X \circ \bd =0$,
$\P_X$ defines a map $\P_X:HFF_*(Y_1,\d_1) \ar HFF_*(Y_2,\d_2)$.
Again if $\a_1 \in H_*(Y_1)$ and $\a_2 \in H_*(Y_2)$ define the same homology class
in $X$, then $\mu(\a_2) \circ \P_X=\P_X\circ \mu(\a_1)$. 

On the other hand, suppose that we have $(Y_1,\d_1,P_1)$ and $(Y_2,\d_2,P_2)$ 
and consider $(Y,\d,P)$ with $Y=Y_1\sqcup Y_2$, $P=P_1\sqcup P_2$ and $\d=\d_1 \sqcup \d_2$. 
One can prove easily that $HFF_*(Y,\d)= HFF_*(Y_1,\d_1) \ox_{\Ct} HFF_*(Y_2,\d_2)$.

Finally, in case that there is a cobordism
between $(Y,\d) \sqcup \, (Y,\d)$ and $(Y,\d)$ we have a map 
\begin{equation}
  HFF_*(Y,\d) \ox_{\Ct} HFF_*(Y,\d) \ar HFF_*(Y,\d),
\label{eqn:f3.2}
\end{equation}
which in some cases it may give an associative and graded commutative ring structure on
$HFF_*(Y,\d)$.
Also note that if there is a cobordism between $(Y,\d) \sqcup \, (Y,\o)$ and $(Y,\d)$
then there will be a map
\begin{equation}
  HF_*(Y) \ox HFF_*(Y,\d) \ar HFF_*(Y,\d),
\label{eqn:f3.3}
\end{equation}
which may lead to a module structure of $HFF_*(Y,\d)$ over $HF_*(Y)$.

\subsection{Relative invariants of $4$-manifolds}
To define relative invariants, let $X$ be a $4$-manifold with
$\bd X=Y$ and $w \in H^2(X;\ZZ)$ such that $w|_Y=w_2\in
H^2(Y;\ZZ/2\ZZ)$. We give $X$ a cylindrical end.
Let $D \subset X$ be a $2$-cycle such that $\bd
D=D\cap Y =\d$ (more accurately, $D \cap (Y \x [0,\infty)) = \d \x
[0,\infty)$). One has the moduli spaces $\cM(X,\r_k)$ and we can
choose generic cycles $V_D^{(i)}$ representing $\mu(D)$ and
intersecting transversely in the top stratum of the
compactification of $\cM(X,\r_k)$ (see~\cite{Fukaya}). Then we have an element
$$
  \p^w(X,D^d)= \sum_{\r_k} (\#
  \cM^{2d}(X, \r_k) \cap V_D^{(1)} \cap \cdots \cap V_D^{(d)}) \r_k 
$$
in $CF_*(Y) \ox H_{2d}(\CP^{\infty}) \subset CFF_*(Y,\d)$.
We remark that this is {\bf not} a cycle. 
Then we set $\p^w(X,D)=\prod_d \p^w(X,D^d)$, which is a
cycle. We also denote by $\p^w(X,D)\in HFF_*(Y,\d)$ 
the Fukaya-Floer homology class
it represents. Alternatively, we denote this same element as 
$$
  \p^w(X,e^{tD}) = \p^w(X,D)=\sum_d \p^w(X,D^d) {t^d \over d!}.
$$
Formally this element lives in $HF_*(Y) \ox \hH_*(\CP^{\infty})$, the $E_3$ term
of the spectral sequence alluded above, but it
represents the same Fukaya-Floer homology class.
The definition of $\p^w(X,D^d)$ depends on some choices~\cite{Fukaya}, but
the homology class $\p^w(X,D)$ only depends on $(X,D)$. Moreover if we
have a homology of $D$ which is the identity on the cylindrical end of
$X$, $\p^w(X,D)$ remains fixed.
Analogously, for any $z \in \AA(X)$, we define $\p^w(X,z\, D^d)
\in CF_*(Y) \ox H_{2d}(\CP^{\infty}) \subset CFF_*(Y,\d)$ and  $\p^w(X,z\, e^{tD})$.
The relevant gluing theorem is~\cite{Fukaya}~\cite{thesis}:

\begin{thm}
\label{thm:f3.fukaya}
  Let $X=X_1 \cup_Y X_2$ and $w \in
  H^2(X;\ZZ)$ with $w|_Y=w_2$. Take 
  $\S \in H^2(X;\ZZ)$ whose Poincar\'e dual lies in the image of
  $H_2(Y;\ZZ) \ar H_2(X;\ZZ)$,
  and satisfies $w \cdot \S \equiv 1\pmod 2$. Put
  $w_i =w|_{X_i} \in H^2(X_i;\ZZ)$. 
  Let $D \in H_2(X)$ decomposed as $D=D_1 +D_2$ with $D_i \subset
  X_i$, $i=1,2$, $2$-cycles with $\bd D_1=\d$, $\bd D_2=-\d$. Choose
  $z_i \in \AA(X_i)$, $i=1,2$. Then
  $$
     \Dws_X(z_1z_2e^{tD})=
     \la\p^{w_1}(X_1,z_1e^{tD_1}),\p^{w_2}(X_2,z_2e^{tD_2})\ra.
  $$
  When $b^+=1$ the invariants are calculated for metrics on $X$ giving a 
  long neck.
\end{thm}

\section{Floer homology of $\Y$}
\label{sec:f4}

We want to specialise to the case relevant to us.  
Let $\S=\S_g$ be a Riemann surface of genus $g \geq 1$ and let $Y=\Y$ be the 
trivial circle bundle over $\S$. Over this $3$-manifold, we fix the $SO(3)$-bundle
with $w_2=\PD[\SS^1]\in H^2(Y;\ZZ/2\ZZ)$, which satisfies the hypothesis of 
section~\ref{sec:f2}. Therefore the instanton Floer homology $HF_*(Y)$ is
well-defined. 
As $Y=\Y$ admits an orientation
reversing self-diffeomorphism, given by conjugation on the $\SS^1$ factor, there
is a Poincar\'e duality isomorphism of $HF^*(Y)$ with the dual of $HF_*(Y)$ 
(this identification will be done systematically and without further notice) 
and a pairing $\la,\ra :HF^*(Y) \ox HF^*(Y) \ar \CC$. 
We introduce a multiplication on $HF^*(Y)$ using the cobordism
between $Y \sqcup \, Y$ and $Y$ given by the $4$-manifold which is
a pair of pants times $\S$. This gives a map 
$HF^*(Y) \ox HF^*(Y) \ar HF^*(Y)$. We shall prove later explicitly that this
is an associative and graded commutative ring structure on $HF^*(Y)$.
As a shorthand notation, we shall write henceforth $HF^*_g=HF^*(Y)$, making
explicit the dependence on the genus $g$ of the Riemann surface $\S$.

The Floer cohomology of $Y=\Y$ has been completely computed thanks to the works
of Dostoglou and Salamon~\cite{DS} and its ring structure has been found by the
author in~\cite{floer} and
turns out to be isomorphic to the quantum cohomology of the moduli space $\M$ of 
stable bundles of odd degree and rank two over $\S$
(with fixed determinant), i.e. $QH^*(\M) \iso HF^*(\S_g \x \SS^1)$, as
the author has proved in~\cite{quantum}.

Here we shall recall the result stated in~\cite{floer}. We fix some notation. 
Let $\{\seq{\g}{1}{2g}\}$ be a symplectic basis of 
$H_1(\S;\ZZ)$ with $\g_i \g_{i+g}=\point$, for $1 \leq i \leq g$. Also 
$x$ will stand for the generator of $H_0(\S;\ZZ)$.
First we recall the usual cohomology ring of $\M$, because of its similarity
with the Floer cohomology $HFF^*_g$ and for later use in section~\ref{sec:f5}.

\subsection{Cohomology ring of $\M$} (see~\cite{King}~\cite{ST}~\cite{quantum})
The ring $H^*(\M)$ is generated by the elements
$$
   \left\{ \begin{array}{l} a= 2\, \mu(\S) \in H^2(\M)
   \\ c_i= \mu (\g_i) \in H^3(\M), \qquad 1 \leq i \leq 2g
   \\ b= - 4 \, \mu(x) \in H^4 (\M)      
    \end{array} \right.
$$
where the map $\mu: H_*(\S) \to H^{4-*}(\M)$ is, as usual, given by
$-{1 \over 4}$ times slanting with the first
Pontrjagin class of the universal $SO(3)$-bundle over $\S\x\M$.
Thus there is a basis $\{f_s\}_{s \in \cS}$ for $H^*(\M)$ with elements of the form
\begin{equation}
   f_s=a^nb^mc_{i_1}\cdots c_{i_r},
\label{eqn:f5.0}
\end{equation}
for a finite set $\cS$ of multi-indices of the
form $s=(n,m; i_1,\ldots,i_r)$, $n,m \geq 0$, $r \geq 0$, $1 \leq
i_1 < \cdots < i_r \leq 2g$. There
is an epimorphism of rings $\AA(\S) \surj H^*(\M)$.
The mapping class group $\Diff(\S)$ acts on $H^*(\M)$, with the action factoring
through the action of the symplectic group $\Spz$ on $\{c_i\}$. 
The invariant part, $H_I^*(\M)$, is generated by
$a$, $b$ and $c=-2 \sum_{i=0}^g c_ic_{i+g}$. Then
\begin{equation}
  \CC[a,b,c]  \surj H^*_I(\M)
\label{eqn:f4.qu4}
\end{equation}
which allows us to write
$$
   H_I^*(\M)= \CC [a,b,c]/I_g,
$$
where $I_g$ is the ideal of relations satisfied by $a$, $b$ and $c$. 
The space $H^3=H^3(\M)$ has a basis $\seq{c}{1}{2g}$, so $\mu:H_1(\S) \isom H^3$.
For $0 \leq k \leq g$, the primitive component of $\L^k H^3$ is
$$
   \L_0^k H^3 = \ker (c^{g-k+1} : \L^k H^3 \ar \L^{2g-k+2} H^3).
$$
The spaces $\L^k_0 H^3$ are irreducible $\Spz$-representations~\cite[theorem 17.5]{Fulton}.
The description of the cohomology ring $H^*(\M)$ is given in the following

\begin{prop}[\cite{ST}~\cite{King}]
\label{prop:f4.homology}
  The cohomology ring of the moduli space $\M$ of stable bundles of odd degree and
  rank two over $\S$ with fixed determinant has a presentation
$$
   H^*(\M)= \bigoplus_{k=0}^{g} \L_0^k H^3 \ox \CC [a, b, c]/I_{g-k},
$$
  where $I_r=(q^1_r,q^2_r,q^3_r)$ and $q^i_r$ are defined recursively by setting
  $q^1_0=1$, $q^2_0=0$, $q^3_0=0$ and then  for all $r \geq 0$ 
  $$
   \left\{ \begin{array}{l} q_{r+1}^1 = a q_r^1 + r^2 q_r^2
   \\ q_{r+1}^2 = b q_r^1 + {2r \over r+1}  q_r^3
   \\ q_{r+1}^3 =  c q_r^1
    \end{array} \right.
  $$
\end{prop}

The basis $\{f_s\}_{s\in\cS}$ of $H^*(\M)$ can be chosen to be as follows. 
Choose, for every $0 \leq k \leq g-1$,
a basis $\{x^{(k)}_i\}_{i \in B_k}$ for $\L^k_0 H^3$. Then
\begin{equation} 
  \{x^{(k)}_i a^n b^m c^r/ k=0,1,\ldots, g-1, \> n+m+r < g-k, \> i \in B_k  \}
\label{eqn:f4.an1}
\end{equation}
is a basis for $H^*(\M)$, as proved in~\cite{ST}.
Also proposition~\ref{prop:f4.homology} gives us the relations for $H^*(\M)$. 
If we set $x^{(k)}_0= c_1 c_2 \cdots c_k \in\L_0^k H^3$, then
the relations are given by 
$$
  x_0^{(k)}q_{g-k}^i, \qquad 1 \leq i\leq 3, \quad 0\leq k \leq g,
$$
and their transforms under the $\Spz$-action.

\subsection{Floer cohomology $HF^*_g$}
The description of the Floer cohomology $HF_g^*=HF^*(Y)$ of $Y=\Y$, where $\S=\S_g$ 
is a Riemann surface of genus $g$, is given in~\cite{floer}. 
Consider the manifold $A=\S \x D^2$, $\S$ times a disc,
with boundary $Y=\Y$, 
and let $\D= \point \x D^2 \subset A$ be the 
horizontal slice. Let $w\in H^2(A;\ZZ)$ be any odd multiple of $\PD [\D]$, 
so that $w|_Y =w_2$. Clearly $\AA(A)= \AA(\S)= \Sym^*(H_0(\S)
\oplus H_2(\S)) \otimes \bigwedge^* H_1(\S)$.
Define the following elements of $HF^*(Y)$ as in~\cite{floer}
\begin{equation}
   \left\{ \begin{array}{l} \a= 2 \, \p^w(A,\S) \in HF^2_g
   \\ \q_i= \p^w(A,\g_i) \in HF^3_g, \qquad 0\leq i \leq 2g
   \\ \b= - 4 \, \p^w(A,x) \in HF^4_g
    \end{array} \right.
\label{eqn:f4.1}
\end{equation}
The relative invariants of section~\ref{sec:f2} give a map
\begin{equation} \begin{array}{ccc}
  \AA(\S) & \ar &  HF^*_g  \\
  z & \mapsto & \p^w(A,z)
\end{array} \label{eqn:f4.2}
\end{equation}
For every $s\in \cS$ and $f_s$ as in~\eqref{eqn:f5.0}, we define
\begin{equation} \begin{array}{lcl}
   z_s&=& \S^n x^m \g_{i_1} \cdots \g_{i_r} \in \AA(\S), \\
   e_s&=& \p^w (A, z_s ) \in HF^*_g.
\end{array} \label{eqn:f4.1.5}
\end{equation}
As a consequence of~\cite[lemma 21]{genusg},
$\{e_s\}_{s\in\cS}$ is a basis for $HF^*_g$. Hence~\eqref{eqn:f4.2} is 
surjective. Now it is easy to check that $\p^w(A,z)\p^w(A,z')=\p^w(A,zz')$, as
for any $s\in\cS$, the gluing theorem~\ref{thm:f2.Floer} implies
$$
  \la\p^w(A,z)\p^w(A,z'),\p^w(A,z_s)\ra =\Dws_{\S\x\CP^1}(zz'z_s)=
   \la\p^w(A,zz'),\p^w(A,z_s)\ra .
$$
In particular this implies that the product of $HF^*_g$ is graded commutative
and associative and that~\eqref{eqn:f4.2} is an epimorphism of rings.
The neutral element of the product is ${\bf 1}=\p^w(A,1)$.
The mapping class group $\Diff(\S)$ acts on $HF^*_g$, with the action factoring
through the action of $\Spz$ on $\{\q_i\}$. It also acts on $\AA(\S)$ 
and~\eqref{eqn:f4.2} is $\Spz$-equivariant.
The invariant part, $(HF_g^*)_I=HF_I^*(Y)$, is generated by
$\a$, $\b$ and $\g=-2 \sum_{i=0}^g \q_i\q_{i+g}$, so that there is an epimorphism
$$
  \Cabg \surj (HF_g^*)_I,
$$
which allows us to write
$$
   (HF_g^*)_I= \Cabg/J_g,
$$
where $J_g$ is the ideal of relations satisfied by $\a$, $\b$ and $\g$. 
As a matter of notation, let $H^3$ denote the $2g$-dimensional vector space generated by
$\seq{\q}{1}{2g}$ in $HF^3$. Then $H^3 \iso H^3(\M)$ and $\p^w(A,\cdot):H_1(\S) \isom
H^3$. No confussion should arise from this multiple use of $H^3$. 
Then from~\cite{floer}, a basis for $HF_g^*$ is given by 
$$ 
  \{x^{(k)}_i\a^a\b^b\g^c/ k=0,1,\ldots, g-1, \> a+b+c < g-k, \> i \in B_k \},
$$
where $x^{(k)}_i \in \L^k_0 H^3$ are interpreted now as Floer products.
The explicit description of $HF^*_g$ is given in~\cite[theorem 16]{floer}

\begin{prop}
\label{prop:f4.floer}
  The Floer cohomology of $Y=\Y$, for $\S=\S_g$ a Riemann surface of genus $g$, and
  $w_2=\PD[\SS^1] \in H^2(Y;\ZZ/2\ZZ)$, has a presentation
  $$
   HF^*(\Y)= \bigoplus_{k=0}^{g} \L_0^k H^3 \ox \Cabg /J_{g-k}.
  $$
  where $J_r=(R^1_r, R^2_r,R^3_r)$ and $R^i_r$ are defined recursively by setting 
  $R^1_0=1$, $R^2_0=0$, $R^3_0=0$ and putting for all $r \geq 0$
  $$
   \left\{ \begin{array}{l} R_{r+1}^1 = \a R_r^1 + r^2 R_r^2
   \\ R_{r+1}^2 = (\b+(-1)^{r+1}8) R_r^1 + {2r \over r+1}  R_r^3
   \\ R_{r+1}^3 = \g  R_r^1
    \end{array} \right.
  $$
\end{prop}

The meaning of this proposition is the following. The Floer (co)homology
$HF_g^*$ is generated as a ring by $\a$, $\b$ and $\q_i$, $1\leq i\leq 2g$,
and the relations are
$$
  x_0^{(k)}R_{g-k}^i, \qquad 1 \leq i\leq 3, \quad 0\leq k \leq g,
$$
where $x_0^{(k)}=\q_1 \q_2 \cdots \q_k \in \L^k_0 H^3$, and the 
$\Spz$-transforms of these. Also if we write
$$
  F_r=\CC[\a,\b,\g]/J_r= (HF^*_r)_I
$$ 
then $HF^*_g=\oplus \L^k_0 H^3 \ox F_{g-k}$.
The quotient $\overline{F}_r={F_r/\g F_r}$ is easily described.

\begin{prop}
\label{prop:f4.fij}
  Let $\overline{F}_r=F_r/\g F_r$, $r\geq 0$.
  Then $\overline{F}_r$ has basis $\a^a\b^b$, $a+b<r$.
  We have $\overline{F}_r =\CC[\a,\b]/ \bar{J}_r$, where $\bar{J}_r=
  (\bar{R}^1_r, \bar{R}^2_r)$, and $\bar{R}^i_r$ are determined by
  $\bar{R}_{0}^1=1$, $\bar{R}_{0}^2=0$ and then recursively for all $r\geq 0$,
$$
   \left\{ \begin{array}{l}\bar{R}_{r+1}^1=\a\bar{R}_r^1 +r^2 \bar{R}_r^2
   \\ \bar{R}_{r+1}^2 = (\b+(-1)^{r+1} 8) \bar{R}_r^1 \end{array}\right.
$$
\end{prop}

\begin{pf}
  The $r+1 \choose 2$ elements $\a^a\b^b$, $a+b<r$, generate 
  $\overline{F}_r$. Also Poincar\'e duality identifies $\overline{F}_r=
  F_r/\g F_r$ with $\ker (\g:F_r \ar F_r)$ which equals $J_{r-1}/J_r$, 
  by~\cite[corollary 18]{floer}. So 
  $$
  \dim\overline{F}_r=\dim (\Cabg/J_r) -\dim (\Cabg/J_{r-1})= 
  \dim F_r-\dim F_{r-1}={r+1 \choose 2},
  $$
  and $\a^a\b^b$, $a+b<r$, form a basis for $\overline{F}_r$. 
\end{pf}

In subsection~\ref{subsec:effective} we shall need the following technical
lemma.

\begin{lem}
\label{lem:f4.redi}
  We have $\bar{J}_r / \bar{J}_{r+1}= 
  \ker (\overline{F}_{r+1}\surj \overline{F}_r) =
  \bigoplus\limits_{-r\leq i\leq r\atop i\equiv r\pmod 2} R_{r+1,i}$, 
  where $R_{r+1,i}$ is a $1$-dimensional vector space such that
  $R_{r+1,i}=\CC[\a,\b]/(\a- 4i\sqrt{-1},\b-8)$ for $r$ even, 
  $R_{r+1,i}=\CC[\a,\b]/(\a- 4i,\b+8)$ for $r$ odd. 
\end{lem}

\begin{pf}
  The first equality follows from the exact sequence
  $$
   {\bar{J}_r \over \bar{J}_{r+1}} \inc \overline{F}_{r+1}=
  {\CC[\a,\b]\over\bar{J}_{r+1}} 
   \surj \overline{F}_r= {\CC[\a,\b] \over \bar{J}_r}.
  $$

  Next we claim that $(\b+(-1)^{r+1} 8) \bar{J}_r \subset \bar{J}_{r+1}
  \subset \bar{J}_r$, $r\geq 0$. The second  inclusion is obvious as 
  $\bar{R}^i_{r+1}$ are written in terms of $\bar{R}^i_r$ by 
  proposition~\ref{prop:f4.fij}.
  The first inclusion follows from $(\b +(-1)^{r+1}8) \bar{R}^1_{r} 
  =\bar{R}^2_{r+1} \in \bar{J}_{r+1}$ and  
  then multiplying the first equation in proposition~\ref{prop:f4.fij}
  by $(\b +(-1)^{r+1} 8)$ to get 
  $(\b +(-1)^{r+1}8) \bar{R}^2_{r} \in \bar{J}_{r+1}$.

  Now $\bar{J}_{r}/\bar{J}_{r+1} = \ker (\b +(-1)^{r+1}8: 
  \overline{F}_{r+1} \ar \overline{F}_{r+1})$. This follows from factoring
  the map $\b +(-1)^{r+1} 8$ as
  $$ 
    \CC[\a,\b] /\bar{J}_{r+1} \surj  \CC[\a,\b] /\bar{J}_{r} 
   \stackrel{\b+(-1)^{r+1} 8}{\hookrightarrow}  \CC[\a,\b] /\bar{J}_{r+1}.
  $$
  The second map is well defined by the claim above and it is 
  a monomorphism since $\a^a\b^b$, $a+b<r$, 
  form a basis for $\CC[\a,\b] /\bar{J}_{r}$,
  and their image under $\b +(-1)^{r+1} 8$ are linearly independent in
  $\CC[\a,\b] /\bar{J}_{r+1}$.
  As $\overline{F}_{r+1}$ is a Poincar\'e duality algebra (being a complete 
  intersection algebra), $\ker (\b +(-1)^{r+1}8:\overline{F}_{r+1} \ar 
  \overline{F}_{r+1})$ is dual to $\overline{F}_{r+1}/(\b +(-1)^{r+1}8)=
  F_{r+1}/(\b +(-1)^{r+1}8,\g)$. Using the computations in the proof 
  of~\cite[proposition 20]{floer}, we get finally 
$$
  \bar{J}_{r}/\bar{J}_{r+1}=F_{r+1}/(\b +(-1)^{r+1}8,\g)= 
  \left\{  \begin{array}{ll}   
   \CC [\a]/\left( (\a^2 +r^2 16) \cdots (\a^2 + 2^2 16)\a \right),
  \qquad & \text{$r$ even} \\
   \CC [\a]/\left( (\a^2 -r^2 16) \cdots (\a^2 - 1^2 16) \right),
  \qquad & \text{$r$ odd}
   \end{array} \right.
$$
  as required.
\end{pf}

\section{Fukaya-Floer homology $HFF_*(\Y,\SS^1)$}
\label{sec:f5}

In this section we are going to describe the Fukaya-Floer (co)homology
of the $3$-manifold $Y=\Y$ with the $SO(3)$-bundle
with $w_2=\PD[\SS^1] \in H^2(Y;\ZZ/2\ZZ)$ and loop 
$\d=\point\x\SS^1 \subset Y=\Y$, together with its ring structure. 
As $Y$ admits an orientation reversing self-diffeomorphism, we can
identify its Fukaya-Floer homology and Fukaya-Floer cohomology through
Poincar\'e duality, which we will do.
From now on we shall fix the genus $g\geq 1$ of $\S$ and 
denote $HFF_g^*=HFF^*(\Y,\SS^1)$.

\subsection{The vector space $HFF_g^*$}
The following argument is taken from~\cite{genusg}.
The spectral sequence computing $HFF^*_g$ has 
$E_3$ term $HF^*_g \ox \Ct$. All the
differentials in this $E_3$ term are of the
form $HF^{\odd}_g \ar HF^{\even}_g$ and
$HF^{\even}_g \ar HF^{\odd}_g$. As $\SS^1$ is invariant under 
the action of the mapping class group
$\Diff (\S)$ on $Y=\Y$, the differentials commute with the
action of $\Diff (\S)$. Since there are elements $f \in
\Diff (\S)$ acting as $-1$ on $H_1(\S)$, we have that $f$
acts as $-1$ on $HF^{\odd}_g$ and as $1$ on
$HF^{\even}_g$. Therefore the differentials are zero.
Analogously for the higher differentials. So the
spectral sequence degenerates in the third term and 
$$
   HFF^*_g= HF^*_g \ox \Ct =HF^*_g[[t]].
$$
The pairing in $HFF^*_g$ is induced from that of $HF^*_g$ by 
coefficient extension to $\Ct$.

For a $4$-manifold $X$ with boundary $\bd X=Y$, 
$w \in H^2(X;\ZZ)$ with $w|_Y=w_2$ and $D \in H_2(X)$ with $\bd D=\SS^1$,
the relative invariants will  be 
$\p^w(X,e^{t D}) \in HF^*_g[[t]]$, i.e. formal power series with 
coefficients in the Floer cohomology $HF^*_g$. 

Recall the manifold $A=\S \x D^2$,
with boundary $Y=\Y$, and let $\D= \point \x D^2 \subset A$ be the 
horizontal slice with $\bd\D=\SS^1$. Let $w
\in H^2(A;\ZZ)$ be any odd multiple of $\PD[\D]$, 
so that $w|_Y =w_2 \in H^2(Y;\ZZ/2\ZZ)$. The elements
\begin{equation}
  \he_s=\p^w(A, z_s\, e^{t\D}) \in HFF_g^*
\label{eqn:f5.nopuse}
\end{equation}
analogous to the elements $e_s$ of~\eqref{eqn:f4.1.5},
for $s \in\cS$, are a basis of $HFF^*_g$ as $\Ct$-module 
(see~\cite[lemma 21]{genusg}).
There is a well defined map $HFF_g^*=HF_g^*\ox\Ct \surj HF_g^*$ formally obtained by equating 
$t=0$. It takes $\p^w(A, z\, e^{t\D}) \mapsto \p^w(A, z)$, for any $z\in\AA(\S)$. 
This map intertwines
the $\mu$ actions on $HFF^*_g$ and $HF^*_g$, and it respects the pairings.

\subsection{The ring $HFF^*_g$}
The ring structure of $HFF^*_g$ comes from the cobordism 
between $(Y,\SS^1) \sqcup (Y,\SS^1)$ and $(Y,\SS^1)$, 
given by the pair of pants
times $(\S,\point)$. This yields
$$
   HFF^*_g \ox HFF^*_g  \ar HFF^*_g 
$$
which is an associative and graded commutative ring structure on $HFF^*_g$.
We prove this as for the case of Floer homology by showing first that 
$\p^w(A, z \, e^{t\D})\p^w(A, z' e^{t\D})=\p^w(A, zz' \, e^{t\D})$, so 
that 
\begin{equation}
\label{eqn:eio}
  \begin{array}{ccl}
   \AA(\S) \ox \Ct &\ar& HFF^*_g \\
   z &\mapsto& \p^w(A, z \, e^{t\D})
  \end{array}
\end{equation}
is a $\Ct$-linear epimorphism of rings. 
The map $HFF^*_g \surj HF^*_g$ mentioned above is a ring epimorphism.
The action of $\mu(\S)$ is Fukaya-Floer multiplication by $\p^w(A,\S \, e^{t\D})$, etc.
We define the following elements of $HFF^*_g$ which are generators as 
$\Ct$-algebra,
\begin{equation}
   \left\{ \begin{array}{l} \ha= 2 \, \q^w(A,\S\,e^{t\D}) \in HFF^2_g
   \\ \hq_i= \q^w(A,\g_i\,e^{t\D}) \in HFF^3_g, \qquad 0\leq i \leq 2g
   \\ \hb= - 4 \, \q^w(A,x\,e^{t\D}) \in HFF^4_g
    \end{array} \right.
\label{eqn:f5.1}
\end{equation}
The mapping class group $\Diff(\S)$ acts on both sides of~\eqref{eqn:eio} 
with the action factoring through an 
action of $\Spz$. The invariants parts surject
\begin{equation}
  \CC[\ha,\hb,\hg]\ox\Ct= \Ct[\ha,\hb,\hg] \surj (HFF_g^*)_I
\label{eqn:f5.2}
\end{equation}
where $\hg =-2 \sum_{i=0}^g \hq_i\hq_{i+g}$. Thus we can write
\begin{equation}
    (HFF^*_g)_I=\Cthabg /\cJ_g,
\label{eqn:f5.3}
\end{equation}
where $\cJ_g$ is the ideal of relations of the generators $\ha$, $\hb$ and $\hg$.
Recall that $t$ has homological degree $2$ and hence cohomological degree
$-2$. The other cohomological degrees are 
$\deg\ha=2$, $\deg\hq_i=3$, $\deg\hb=4$ and $\deg\hg=6$. 

The ring structure of $HFF^*_g$, which is in some sense equivalent to the
determination of the kernel of~\eqref{eqn:f5.2}, runs closely parallel to
the arguments in~\cite{floer} to find out the ring structure of 
$HF^*_g=HF^*(Y)$. We recommend the reader have~\cite{floer} at hand.

Consider the ring $H^*(\M)[[t]]$, where $t$ is given degree $-2$. 
The elements in $H^i(\M)[[t]]$ are thus sums 
$\sum_{n\geq 0} s_{i+2n} t^n$, where $\deg(s_{i+2n})=
i+2n$. Note that all such elements are finite sums, although $H^i(\M)[[t]]\neq 0$ for 
arbitrarily negative $i$. The analogue of~\cite[theorem 5]{floer} is

\begin{prop}
\label{prop:f5.deform}
  Denote by  $*$ the product induced in $H^*(\M)[[t]]$ by the product in
  $HFF^*_g$ 
  under the $\Ct$-linear isomorphism $H^*(\M)[[t]] \isom HFF^*_g$ given by
  $f_s\mapsto \he_s$, 
  $s \in \cS$. Then $*$ is a deformation of the cup-product
  graded modulo $4$, i.e. for $f_1 \in H^i (\M)[[t]]$, 
  $f_2 \in H^j(\M)[[t]]$, it is
  $f_1 * f_2 = \sum_{r \geq 0} \P_r(f_1,f_2)$, where $\P_r \in
  H^{i+j-4r}(\M)[[t]]$ and 
  $\P_0= f_1 \cup f_2$.
\end{prop}

\begin{pf}
  To start with, let us fix some notation. The choice of basis~\eqref{eqn:f5.nopuse}
  gives a splitting $\imath:
  H^*(\M) \ar \AA(\S)$, $f_s \mapsto z_s$, satisfying the property that
  $f \mapsto \p^w(A, \imath (f)\, e^{t\D})$ is the isomorphism of the statement.

  Now we claim that for any $s, s' \in \cS$ it is
  \begin{equation}
   \la\he_s, \he_{s'}\ra  = \Dws_{\S \x \CP^1}(z_sz_{s'}e^{t\CP^1}) = -\la f_s,
  f_{s'}\ra +  O(t^{(6g-6-(\deg(f_s) + \deg(f_{s'})))/2+1}),
  \label{eqn:f5.claim}
  \end{equation}
  where $O(t^r)$ means any element in $t^r\Ct$ (note that~\eqref{eqn:f5.claim} vanishes when
  $\deg(f_s) + \deg(f_{s'}) \not\equiv 0 \pmod 2$).
  If $\deg(f_s) + \deg(f_{s'}) > 6g-6$ then the statement is vacuous. 
  For $\deg(f_s) + \deg(f_{s'})  \leq 6g-6$ it follows from the fact that 
  the dimensions of the moduli spaces of anti-self-dual connections 
  on $\S \x
  \CP^1$ are $6g-6+4r$, $r \geq 0$, and the $(6g-6)$-dimensional moduli space is
  $\M$, as remarked in~\cite{floer}, so that for $\deg(f_s) + \deg(f_{s'}) +2d=6g-6$, it is
  $\Dws_{\S \x \CP^1}(z_sz_{s'}(\CP^1)^{d})=0$, unless $d=0$, and in 
  that case it gives $-\la f_s, f_{s'}\ra $ (the
  minus sign is due to the different convention orientation for 
  Donaldson invariants). 

  We shall check the statement of the proposition on basic elements
  $f_s$ and $f_{s'}$ of degrees $i$ and $j$ respectively.
  Put $f_s * f_{s'}= \sum_{m \leq M} g_m$, where $g_m \in H^m(\M)[[t]]$
  and $g_M \neq 0$
  is the leading term.
  By definition, $\he_s\he_{s'}= \sum_{m \leq M} \hat g_m$
  (with $\hat g_m \in HFF^*_g$ corresponding
  to $g_m$ under the isomorphism of the statement).

  Suppose $M > i+j$. Then let $f t^r$, $f\in H^*(\M)$, be the non-zero
  monomial in $g_M$ with
  minimum $r$. So $f$ has degree $M+2r$. Pick $f' \in H^{6g-6-(M+2r)}(\M)$
  with $\la f,f'\ra=-1$ in $H^*(\M)$. Let $z,z' \in \AA(\S)$ be the elements 
  corresponding to
  $f , f' \in H^*(\M)$ under the splitting $\imath$. Then by~\eqref{eqn:f5.claim}
  $$
    \la t^r\p^w(A,ze^{t\D}), \p^w(A, z'e^{t\D})\ra  = t^r + O(t^{r+1}),
  $$
  so $\la\hat g_M, \p^w(A, z'e^{t\D})\ra =t^r +O(t^{r+1})$. For $m <M$, it must be
  $\la\hat g_m, \p^w(A, z'e^{t\D})\ra =O(t^{r+1})$ by~\eqref{eqn:f5.claim} again, so finally
  $$
    \la\he_s\he_{s'}, \p^w(A, z'e^{t\D})\ra  = t^r + O(t^{r+1}).
  $$
  On the other hand, as $\deg(f_s)+\deg(f_{s'})+\deg(f') <6g-6-2r$, it is
  $$
    \la\he_s\he_{s'}, \p^w(A, z'e^{t\D})\ra =\Dws_{\S \x
  \CP^1}(z_sz_{s'}z'e^{t\CP^1})=O(t^{r+1}),
  $$
  which is a contradiction. It must be $M \leq i+j$.

  For $m=i+j$, put $g_m=G_{i+j} +tG_{i+j+2} + \cdots$, where $G_{i+j+2r} \in H^{i+j+2r}(\M)$.
  Pick any $f_{s''}$ of degree $6g-6-m$. Clearly $\Dws_{\S \x
   \CP^1}(z_sz_{s'}z_{s''}e^{t\CP^1})=-\la f_s f_{s'}, f_{s''}\ra +O(t)$. Also 
$$
  \Dws_{\S\x \CP^1}(z_sz_{s'}z_{s''}e^{t\CP^1})=
  \la\he_{s}\he_{s'},\he_{s''}\ra 
  =\la\hat g_m, \he_{s''}\ra +O(t)=-\la g_m, f_{s''}\ra +O(t).
$$
  So $\la G_{i+j},f_{s''}\ra =\la f_sf_{s'},f_{s''}\ra $, for arbitrary $f_{s''}$, and hence
  $G_{i+j}=f_sf_{s'}$. 

  To check that $G_{i+j+2r}=0$ for $r>0$, 
  pick any $f_{s''}$ of degree $6g-6-(m+2r)$. By~\eqref{eqn:f5.claim} it is
  $\Dws_{\S \x \CP^1}(z_sz_{s'}z_{s''}e^{t\CP^1})=O(t^{r+1})$ and
$$
  \Dws_{\S \x \CP^1}(z_sz_{s'}z_{s''}e^{t\CP^1})=\la\he_{s}\he_{s'},\he_{s''}\ra =
  \la\hat g_m, \he_{s''}\ra +O(t^{r+1})=- \la G_{i+j+2r}, f_{s''}\ra  t^r+O(t^{r+1}),
$$ 
  So $\la G_{i+j+2r},f_{s''}\ra =0$, i.e. $G_{i+j+2r}=0$. 
\end{pf}

From~\eqref{eqn:f4.an1} a basis of $H^*(\M)[[t]]$ as $\Ct$-module is given by 
$$
  \{x^{(k)}_ia^nb^mc^r/ k=0,1,\ldots, g-1, \> n+m+r < g-k, \> i \in B_k \}.
$$
Recalling $x^{(k)}_0=c_1c_2 \cdots c_k \in\L_0^k H^3$,  
a complete set of relations satisfied in $H^*(\M)$
are $x_0^{(k)} q^i_{g-k}$, $i=1,2,3$, $0 \leq k \leq g$, and the $\Spz$-transforms of these.
Now identifying $H^3$ with the $2g$-dimensional subspace of $HFF^3_g$ generated
by $\seq{\hq}{1}{2g}$, proposition~\ref{prop:f5.deform} implies that the set
$$ 
  \{x^{(k)}_i\ha^a\hb^b\hg^c/ k=0,1,\ldots, g-1, \> a+b+c < g-k, \> i \in B_k \},
$$
where $x^{(k)}_i \in \L^k_0 H^3 \subset HFF_g^*$, 
is a basis
for $HFF_g^*$ as $\Ct$-module, where Fukaya-Floer multiplication is understood. Also from
proposition~\ref{prop:f4.homology}, we can write
$$
  H^*(\M)[[t]]=\bigoplus_{k=0}^g \L^k_0 H^3 \ox {\Ct[a,b,c] \over 
  (q_{g-k}^1,q_{g-k}^2,q_{g-k}^3)}.
$$
The products in both $H^*(\M)[[t]]$ and $HFF_g^*$ are $\Spz$-equivariant and the 
isomorphism in the statement of proposition~\ref{prop:f5.deform} is also
$\Spz$-equivariant. Then 
we can use the arguments in the proof of~\cite[proposition 16]{quantum}
to write
$$
  HFF_g^*= \bigoplus_{k=0}^g \L^k_0 H^3 \ox {\Cthabg \over 
 (\cR_{g-k}^1,\cR_{g-k}^2,\cR_{g-k}^3)},
$$
where if we put $x^{(k)}_0=\hq_1\hq_2 \cdots \hq_k \in\L_0^k H^3$, then
$x^{(k)}_0 \cR_{g-k}^i$, $i=1,2,3$, $0 \leq k \leq g$, and their $\Spz$-transforms,
are a complete set of relations for $HFF_g^*$.
More explicitly, we decompose $HFF_g^*=\bigoplus\limits_{k=0}^g \cV_k$, 
where $\cV_k=\L^k_0 H^3 \ox \cF_{g-k}$ is the image of
$$
  \L^k_0 H^3 \ox \Cthabg \ar HFF_g^*,
$$
so in particular, the invariant part is $\cV_0=(HFF_g^*)_I$ and $\cV_g=0$. Then 
$$
  HFF_g^*= \bigoplus_{k=0}^g \L^k_0 H^3 \ox \cF_{g-k},
$$
where
\begin{equation}
   \cF_{g-k}={\Cthabg \over \cJ_{g-k}},
\label{eqn:cFr}
\end{equation}
for $0\leq k\leq g$, where the generators of the ideal 
$\cJ_{g-k} \subset \Cthabg$ are obtained by writing 
$q^1_{g-k}$, $q^2_{g-k}$, $q^3_{g-k}$ in terms of the Fukaya-Floer product 
(see~\cite{ST2}
for an analogous argument in the study of quantum cohomology)
$$
    q^i_{g-k} = \sum  c^i_{abcd} \ha^a\hb^b\hg^c t^d,
$$
where the sum runs for $a+b+c<g-k$, $d\geq 0$, $2a+4b+6c-2d =\deg q^i_{g-k}-4r$, $r>0$
and
$c^i_{abcd} \in \CC$. So $\cJ_{g-k}= (\cR_{g-k}^1,\cR_{g-k}^2,\cR_{g-k}^3)$ with
$$
   \cR^i_{g-k}= q^i_{g-k} -\sum c^i_{abcd} \ha^a\hb^b\hg^c t^d.
$$
The elements $\cR^i_{g-k}$ are uniquely defined by the folowing two conditions: 
$x_0^{(k)}\cR^i_{g-k}=0 \in HFF^*_g$ and $\cR^1_{g-k}$ (respectively $\cR^2_{g-k}$,
$\cR^3_{g-k}$) equals $\ha^{g-k}$ (respectively $\ha^{g-k-1}\hb$, $\ha^{g-k-1}\hg$) plus terms
of the form $\ha^a\hb^b\hg^c t^d$ with $a+b+c <g-k$. 
Note that they might depend, in principle, not only on $g-k$ but also on the genus $g$ 
(which was fixed throughout this section). In particular $\cR^1_0=1$,
$\cR^2_0=0$ and $\cR^3_0=0$.

\begin{lem}
\label{lem:f5.pr}
  For $0 \leq r \leq g-1$, it is $\hg \cJ_r \subset \cJ_{r+1}\subset \cJ_r$.
\end{lem}

\begin{pf}
Analogous to~\cite[lemma 17]{quantum}.
\end{pf}

\begin{thm}
\label{thm:f5.main}
  Fix $g \geq 1$. Let $\S=\S_g$ be a Riemann surface of genus $g$. The
  Fukaya-Floer cohomology $HFF_g^*=HFF^*(\Y,\SS^1)$
  has a presentation
  $$
   HFF^*_g= \bigoplus_{k=0}^g \L_0^k H^3 \ox \Cthabg /\cJ_{g-k}.
  $$
  where $\cJ_r=(\cR^1_r, \cR^2_r,\cR^3_r)$ and $\cR^i_r$ are defined 
  recursively by setting 
  $\cR^1_0=1$, $\cR^2_0=0$, $\cR^3_0=0$ and putting, for all $0\leq
  r \leq g-1$,
\begin{equation}
   \left\{ \begin{array}{l} \cR_{r+1}^1 = (\ha+f_{11}(t)) \cR_r^1 + r^2 
     (1+f_{12}(t)) \cR_r^2 +f_{13}(t)\cR_r^3
   \\ \cR_{r+1}^2 = (\hb+(-1)^{r+1} 8+f_{21}(t)) \cR_r^1 +f_{22}(t)\cR_r^2+ 
   ({2r \over r+1}+f_{23}(t)) \cR_r^3
    \\ \cR_{r+1}^3 = \hg  \cR_r^1 
    \end{array} \right.
\label{eqn:f5.tex}
\end{equation}
  for some (unknown) functions $f_{ij}^{r,g}(t) \in t\Cthabg$, 
  dependent on $r$ and $g$. Moreover $f_{ij}$ are such that 
  $f_{11}\cR^1_r+f_{12}\cR^2_r+f_{13}\cR^3_r$ and 
  $f_{21}\cR^1_r+f_{22}\cR^2_r+f_{23}\cR^3_r$ 
  are both $\Ct$-linear combinations of the monomials
  $\ha^a\hb^b\hg^c$, $a+b+c< r+1$. 
\end{thm}

\begin{pf}
We only have to prove the recurrence stated in~\eqref{eqn:f5.tex}, 
which is similar to~\cite[theorem 10]{floer}.
The inclusion $\hg \cJ_{r} \subset \cJ_{r+1}$ says that
$\hg \cR_r^1$ must be in $\cJ_{r+1}$, so it must coincide with $\cR_{r+1}^3$.
Now the inclusion $\cJ_{r+1}\subset \cJ_r$ 
implies the recurrence as written in~\eqref{eqn:f5.tex} with 
$f_{ij} \in \Cthabg$. Lastly, the $\Spz$-equivariant epimorphism 
$HFF_g^* \surj HF_g^*$ yields that $\cR^i_g$ reduces to $R^i_g$ when 
we set $t=0$. Thus the functions $f_{ij}$ are multiples of $t$.

The last sentence follows from
the fact that $\cR_{r+1}^1$ is written as a series with leading term
$\ha^{r+1}$ plus terms of the form $\ha^a\hb^b\hg^c t^d$, $a+b+c< r+1$, 
and that $\cR_{r+1}^2$ is written as a series with leading term
$\ha^{r}\hb$ plus terms of the form $\ha^a\hb^b\hg^c t^d$, $a+b+c< r+1$.
\end{pf}

We give the following two results, whose proofs are left to the reader, 
for completeness.

\begin{cor}
\label{cor:f5.nSS1}
  Fix $g \geq 1$. Let $\S=\S_g$ be a Riemann surface of genus $g$. Let $n \in\ZZ$. The
  Fukaya-Floer cohomology $HFF^*(\Y,n \SS^1)$
  has a presentation
  $$
   HFF^*(\Y,n \SS^1)= \bigoplus_{k=0}^g \L_0^k H^3 \ox \Cthabg /\cJ_{g-k}.
  $$
  where $\cJ_r=(\cR^1_r, \cR^2_r,\cR^3_r)$ and $\cR^i_r$ are defined 
  recursively by setting 
  $\cR^1_0=1$, $\cR^2_0=0$, $\cR^3_0=0$ and putting, 
  for all $0 \leq r\leq g-1$,
$$
   \left\{ \begin{array}{l} \cR_{r+1}^1 = (\ha+f_{11}(nt)) 
   \cR_r^1 + r^2 
     (1+f_{12}(nt)) \cR_r^2 +f_{13}(nt)\cR_r^3
   \\ \cR_{r+1}^2 = (\hb+(-1)^{r+1} 8+f_{21}(nt)) \cR_r^1 +f_{22}(nt)\cR_r^2+ 
   ({2r \over r+1}+f_{23}(nt)) \cR_r^3
    \\ \cR_{r+1}^3 = \hg  \cR_r^1 
    \end{array} \right.
$$
\end{cor}

In particular, for $n=0$ we recuperate proposition~\ref{prop:f4.floer}.

\begin{prop}
\label{prop:f5.reader}
  Let $n \in \ZZ$ and consider $HFF^*(\Y,n \SS^1)$. Then for any $\a \in 
  H_0(\S)$ or $\a\in H_1(\S)$, the 
  action of $\mu (\a \x\SS^1)$ in $HFF^*(\Y,n \SS^1)$ is zero.
\end{prop}

\subsection{Reduced Fukaya-Floer homology}
\label{subsec:rHFF}
We introduce the following space associated to $HFF_g$ in order to 
understand the gluing theory of
$4$-manifolds with $b^+>1$ which are of strong simple type.

\begin{defn}
\label{def:rHFF}
  We define the {\em reduced Fukaya-Floer homology of $\Y$\/} to be
  $$
     \overline{HFF}{}^*_g= HFF^*_g / (\hb^2-64, \seq{\hq}{1}{2g})
  \iso \ker(\hb^2-64) \cap \bigcap_{1\leq i\leq 2g} \ker \hq_i,
  $$
  where the last isomorphism is Poincar\'e duality. Note that
  $\overline{HFF}{}^*_g= (HFF^*_g)_I / (\hb^2-64,\hg)$.
\end{defn}

Suppose that $X_1$ is a $4$-manifold with boundary $\bd X_1=Y$ 
(and $w\in H^2(X;\ZZ)$ satisfies $w|_Y=w_2=\PD [\SS^1]$)
such that $X=X_1 \cup_Y A$ is a closed $4$-manifold with $b^+>1$
and of strong simple type. Then 
$$
  \p^w(X_1, z_1 e^{tD_1}) \in \ker(\hb^2-64) \cap 
  \bigcap_{1\leq i\leq 2g} \ker \hq_i \iso \overline{HFF}{}^*_g,
$$
for any $z_1 \in \AA(X_1)$ and any $D_1\subset X_1$ with $\bd D_1=\SS^1$. 
Indeed $\la(\hb^2-64)\p^w(X_1, z_1e^{tD_1}), \he_s\ra = \Dws_X(16 z_1z_s
(x^2-4) e^{tD})=0$, for any $s\in\cS$, where $D=D_1+\D \in H_2(X)$. Then
$(\hb^2-64)\p^w(X_1, z_1e^{tD_1})=0$. Analogously 
$\hq_i \p^w(X_1, z_1e^{tD_1})=0$, for $1\leq i \leq 2g$.

\begin{prop}
\label{prop:f5.10}
  For $0\leq k \leq g-1$, there exists a non-zero vector $v \in (HFF^*_g)_I$ 
  such that
  \begin{eqnarray*}
    \ha v &=& \left\{ \begin{array}{ll} 
    (\pm 4(g-k-1)+2t) \, v & \text{$g-k$ even} \\
    (\pm 4(g-k-1) \sqrt{-1}-2t) \, v \qquad & \text{$g-k$ odd} 
    \end{array} \right. \\
    \hb v &=& (-1)^{g-k-1} 8 \, v \\
    \hg v & =& 0
  \end{eqnarray*}  
\end{prop}

\begin{pf}
  This is an extension of~\cite[proposition 12]{floer}.
  We shall construct the vector corresponding to the plus sign, the other
  being analogous. We have the following cases:
\begin{itemize}
  \item $0=k<g-1$. We shall look for $v \in (HFF_g^*)_I$ constructed
  as the relative invariants of a particular $4$-manifold (see section~\ref{sec:f3}).
  For finding such a vector $v$ we use the same manifold as in the proof
  of~\cite[proposition 12]{floer}. This is a
  $4$-manifold $X=C_g$ with an embedded Riemann surface $\S$ of genus $g$ and
  self-intersection zero, and $w \in H^2(X;\ZZ)$ with $w\cdot \S \equiv 1 \pmod 2$.
  Such $X$ is of simple type, with $b_1=0$, $b^+>1$. 
  Suppose for simplicity that $g-k=g$ is even (the other case is analogous). Then 
  the Donaldson invariants of $X$ are
\begin{equation}
  \Dws_{X}(e^{\a})= -2^{3g-5} e^{Q(\a)/2} e^{K \cdot \a} +
   (-1)^g 2^{3g-5} e^{Q(\a)/2} e^{-K \cdot \a},
\label{eqn:f5.nose}
\end{equation} 
  for a single basic class $K\in H^2(X;\ZZ)$ with $K\cdot \S=2g-2$.
  Let $X_1$ be $X$ with a small open tubular neighbourhood of $\S$ removed, so that
  $X=X_1 \cup_Y A$. Consider $D\subset X$ intersecting transversely $\S$ in just one
  positive point. Let $D_1=X_1 \cap D\subset X_1$, so that $\bd D_1=\SS^1$ and 
  $D=D_1 +\D$. Then
  set 
  $$
  v = \p^w(X_1,(\S + 2g-2-t)e^{tD_1}) \in HFF^*(\Y,\SS^1)_I.
  $$
  Let us prove that this $v$ does the job. For any $z_s=\S^n x^m 
  \g_{i_1} \cdots \g_{i_r}$, we compute from~\eqref{eqn:f5.nose} that
  $\la v, \he_s\ra =\la \p^w(X_1,(\S + 2g-2-t)e^{tD_1}), \p^w(A, z_s e^{t\D})\ra 
  =  \Dws_{X}((\S+2g-2-t) z_s e^{tD})$ is equal to
$$
   \left\{ \begin{array}{ll} 0, & r >0 \\ -2^{3g-4}(2g-2)(2g-2+t)^n 2^m\, e^{Q(tD)/2+tK\cdot D},
   \qquad & r=0 \end{array} \right.
$$
  Then $\la \ha v , \he_s\ra  = \la \p^w(X_1,2 \,\S(\S + 2g-2-t)e^{tD_1}), \p^w(A, z_s e^{t\D})\ra = 
  \Dws_X((\S+2g-2-t) 2 \,\S\, z_s e^{t D})= (4g-4+2t) \la v, \he_s\ra $, for all $s \in 
  \cS$. Then $\ha v = (4g-4+2t) v$. Analogously,
  $\hg v =0$ and $\hb v= -8 v$.

  \item $0<k<g-1$. The same argument above for genus $g-k$ produces
  a $4$-manifold $C_{g-k}$ with an embedded Riemann surface $\S_{g-k}$ of genus
  $g-k$ and self-intersection zero with a single basic class $K\in H^2(X;\ZZ)$ with 
  $K\cdot \S_{g-k}=2(g-k)-2$. 
  Let now $X=C_{g-k} \# k\, \SS^1 \x \SS^3$ (performing the connected sum well appart from
  $D$). Consider the torus $\SS^1\x\SS^1 \subset \SS^1\x
  \SS^3$ and the internal connected sum $\S=\S_g =\S_{g-k} \# k\, \SS^1\x\SS^1 \subset
  X=C_{g-k} \# k\, \SS^1\x\SS^3$. When choosing the basis of $H_1(\S;\ZZ)$, we arrange 
  $\g_1,\ldots, \g_k$ such that $\g_i=\SS^1 \x \point$ in the $i$-th copy
  $\SS^1\x\SS^3$. Suppose for instance that $g-k$ is even. Then by lemma~\ref{lem:f5.S1xS3}
  below it is, for any $\a \in H_2(C_{g-k})=H_2(X)$,
$$
  \Dws_{X}(\g_1\cdots\g_k e^{\a})= c^k \,\Dws_{C_{g-k}}(e^{\a}) =
   c^k \left( -2^{3(g-k)-5} e^{Q(\a)/2} e^{K \cdot \a} +
   (-1)^{g-k} 2^{3(g-k)-5} e^{Q(\a)/2} e^{-K \cdot \a} \right),
$$
  with $w \in H^2(C_{g-k};\ZZ)$ as in the first case. 
  Write again
  $X=X_1 \cup_Y A$ and consider $D\subset X$ intersecting transversely $\S$ in one 
  point with $D\cdot \S=1$, so that 
  $D=D_1 +\D$ with $\bd D_1=\SS^1$. Then the element
$$
   v=\p^w(X_1, (\S+2(g-k)-2-t)\g_1\cdots \g_k e^{tD_1})\in (HFF^*_g)_I
$$
  satisfies the required properties. Note that $v$ is invariant since has
  only non-zero pairing with elements in $(HFF_g^*)_I \subset HFF_g^*$.

  \item $k=0$ and $g=1$. Let $S$ be the $K3$ surface and fix
  an elliptic fibration for $S$, whose generic fibre is an embedded torus $\S=\TT^2$. The 
  Donaldson invariants are, for $w \in H^2(S; \ZZ)$ with $w \cdot \S \equiv 1
  \pmod 2$, $\Dws_S (e^{tD})=-e^{-Q(tD)/2}$.
  Fix $D \subset S$ cutting $\S$ transversely in one point such that $\S\cdot D=1$.
  Then $\Dws_S (\S \, e^{tD})= te^{-Q(tD)/2}$.
  Let $S_1$ be the complement of a small open tubular neighbourhood of $\S$ in $S$
  and $D_1=S_1\cap D \subset S_1$, so that $\bd D_1=\SS^1$. 
  Then $v=\p^w(S_1, e^{tD_1})$ generates $HFF_1^*$ and $\p^w(S_1, \S \, e^{tD_1})=
  -t \p^w(S_1, e^{tD_1})$, so that $\ha v =-2t v$. Analogously
  $\hb v=8 v$ and $\hg v=0$.

  \item $0 <k=g-1$. We use the same trick as in the second case, considering
  the $K3$ surface connected sum with $k$ copies of $\SS^1\x\SS^3$.
\end{itemize}
\end{pf}

\begin{lem}
\label{lem:f5.S1xS3}
   Let $X$ be a $4$-manifold with $b^+> 1$, and $z \in \AA(X)$. Then consider
   $\tilde X=X\# \SS^1\x\SS^3$ and $\g=\SS^1\x\point \subset \SS^1\x\SS^3$ the natural
   generator of the fundamental group of $\SS^1\x\SS^3$. We can view $\g$ as
   an element of $\AA(\tilde X)$. For any $w \in H^2(X;\ZZ)=H^2(\tilde X;\ZZ)$,
   it is $D^w_{\tilde X}(\g \, z)=c \, D^w_X (z)$, where $c$ is a universal constant. 
\end{lem}

\begin{pf}
  Consider the moduli space $\cM^{w,\kappa}_{\tilde X}$ of ASD connections over $\tilde X$
  of dimension $d=3+\deg(z)$, $\kappa$ denoting the charge~\cite{KM}.
  Then there is a choice of generic cycles $V_z$, $V_{\g}$ in
  $\cM^{w,\kappa}_{\tilde X}$ such that
  $$
    M= \cM^{w,\kappa}_{\tilde X} \cap V_z
  $$
  is smooth $3$-dimensional and compact, and $D^w_{\tilde X}(\g \, z)=\# (M\cap V_{\g})$.
  For metrics giving a long neck to the connected sum $\tilde X=X\# \SS^1\x\SS^3$, the
  usual counting arguments give that the only possible distribution of charges of 
  limitting connections are $\kappa$ on the $X$ side and $0$ on the $\SS^1\x\SS^3$.
  Now recall that the moduli space of flat $SO(3)$-connections on $\SS^1\x\SS^3$ is
  $\Hom(\pi_1(\SS^1\x\SS^3), SO(3))=SO(3)$, hence
  $$
  M= (\cM^{w,\kappa}_{X}) \cap V_z \x SO(3),
  $$
  where $\cM^{w,\kappa}_{X} \cap V_z$ consists of $D^w_X (z)$ points (counted with signs).
  The description of the cycle $V_{\g}$ given in~\cite{KM} implies that there is a universal
  constant $c=\# SO(3)\cap V_{\g}$ yielding the statement of the lemma. 
\end{pf}

\begin{thm}
\label{thm:f5.rHFF}
  $\overline{HFF}{}^*_g$ is a free $\Ct$-module of rank $2g-1$. Moreover
  $\overline{HFF}{}^*_g=\bigoplus\limits_{i=-(g-1)}^{g-1} R_{g,i}$, where
  $R_{g,i}$ are free $\Ct$-modules of rank $1$ such that 
  for $i$ odd, $\ha=4i +2t$ and $\hb=-8$ in $R_{g,i}$. For $i$ even, 
  $\ha=4i\sqrt{-1} -2t$ and $\hb=8$ in $R_{g,i}$. 
\end{thm}

\begin{pf}
  Since $\ha^a\hb^b\hg^c$, $a+b+c<g$ form a basis for $(HFF^*_g)_I$, 
  we have that $\ha^a\hb^b$, $a+b<g$, $b=0,1$, generate 
  $\overline{HFF}{}^*_g$ as $\Ct$-module. Therefore the rank of
  $\overline{HFF}{}^*_g$ is less or equal than $2g-1$.
  Now using Poincar\'e duality in $HFF^*_g$, we have that the dual of
  $\overline{HFF}{}^*_g$ is 
  \begin{equation}
   \ker (\hb^2-64) \cap \ker \hq_1 \cap \cdots \cap \ker \hq_{2g} \subset
   HFF^*_g.
  \label{eqn:f5.k}
  \end{equation}
  By proposition~\ref{prop:f5.10} there are at least $2g-1$ independent
  vectors in~\eqref{eqn:f5.k}, so the rank of $\overline{HFF}{}^*_g$ is
  exactly $2g-1$ and it must be a free $\Ct$-module. The $2g-1$ 
  eigenvalues of $(\ha,\hb)$ given by proposition~\ref{prop:f5.10} 
  provide the decomposition in the statement.
\end{pf}

\subsection{Effective Fukaya-Floer homology} 
\label{subsec:effective}
We introduce the following
space associated to $HFF_g^*$ in order to understand the gluing theory of
general $4$-manifolds which have $b^+>1$.

\begin{defn}
\label{def:f5.eHFF}
  We define the effective Fukaya-Floer homology as the sub-$\Ct$-module
  $\whff \subset HFF^*_g$ generated by all $\p^w(X_1, z_1 e^{tD_1})$, for
  all $4$-manifolds $X_1$ with boundary $\bd X_1=Y$ such that $X=X_1\cup_Y
  A$ has $b^+>1$, $z_1 \in \AA(X_1)$, $D_1 \subset X_1$ with $\bd D_1=\SS^1$
  and $w\in H^2(X_1;\ZZ)$ with $w|_Y=w_2=\PD[\SS^1]$.
\end{defn}

The action of $\Spz$ on $HFF_g^*$ restricts to an action on
$\whff$. Also $\ha$, $\hb$, $\seq{\hq}{1}{2g}$ (and hence 
$\hg$) act on $\whff$ by multiplication.
Now 
 $$
  \whff=\bigoplus_{k=0}^g \L^0_k H^3 \ox \widetilde{\cF}_{g-k} \subset
  \bigoplus_{k=0}^g \L^k_0 H^3 \ox \cF_{g-k},
 $$
with $\widetilde{\cF}_r \subset \cF_r$, $0\leq r\leq g$. We aim to find
the eigenvalues of $(\ha,\hb,\hg)$ in $\widetilde{\cF}_r$. Some preliminary
information on $\cF_r$ is needed. The following result
is proved analogously to proposition~\ref{prop:f4.fij},

\begin{prop}
\label{prop:f5.redi}
  Let $\overline{\cF}_r=\cF_r/\hg \cF_r$, $0\leq r\leq g$.
  Then $\overline{\cF}_r$ is a free $\Ct$-module with basis
  $\ha^a\hb^b$, $a+b<r$.
  We have $\overline{\cF}_r =\Ct[\ha,\hb]/ \bar{\cJ}_r$, where $\bar{\cJ}_r=
  (\bar{\cR}^1_r, \bar{\cR}^2_r)$, and $\bar{\cR}^i_r$ are determined by
  $\bar{\cR}_{0}^1=1$, $\bar{\cR}_{0}^2=0$ and for $0 \leq r \leq g-1$, 
$$
   \left\{ \begin{array}{l} \bar{\cR}_{r+1}^1 = (\ha+\bar{f}_{11}) 
     \bar{\cR}_r^1 + r^2 (1+\bar{f}_{12}) \bar{\cR}_r^2
   \\ \bar{\cR}_{r+1}^2 = (\hb+(-1)^{r+1} 8+\bar{f}_{21}) \bar{\cR}_r^1 +
   \bar{f}_{22}\bar{\cR}_r^2 \end{array}\right.
$$
  for some functions $\bar{f}_{ij}(t)\in t\Ct[\ha,\hb]$.
  Moreover $\bar f_{ij}$ are such that 
  $\bar f_{11}\bar\cR^1_r+\bar f_{12}\bar\cR^2_r$ and 
  $\bar f_{21}\bar\cR^1_r+\bar f_{22}\bar\cR^2_r$ 
  are both $\Ct$-linear combinations of the monomials
  $\ha^a\hb^b$, $a+b< r+1$. 
\hfill $\Box$
\end{prop}

\begin{lem}
\label{lem:f5.redi}
  We have $\bar{\cJ}_r / \bar{\cJ}_{r+1}= 
  \ker (\overline{\cF}_{r+1}\surj \overline{\cF}_r) =
  \bigoplus\limits_{-r\leq i\leq r\atop i\equiv r\pmod 2} R_{r+1,i}$, 
  where $R_{r+1,i}$ is a free $\Ct$-module of rank $1$. For $r$ even,
  $\ha=4i\sqrt{-1}+O(t)$ and $\hb=8+O(t)$ in $R_{r+1,i}$. For $r$ odd,
  $\ha=4i+O(t)$, $\hb=-8+O(t)$ in $R_{r+1,i}$.
\end{lem}

\begin{pf}
  The natural map $HFF^*_g \surj HF^*_g$ given by equating $t=0$ 
  and lemma~\ref{lem:f4.redi} give the following 
  commutative diagram with exact rows
\begin{equation*}
\begin{array}{ccccc}
   {\bar{\cJ}_r / \bar{\cJ}_{r+1}} &\inc& \overline{\cF}_{r+1} &\surj &\overline{\cF}_r \\
        \downarrow & & \downarrow & &\downarrow \\
   {\bar{J}_r / \bar{J}_{r+1}} &\inc &\overline{F}_{r+1}& \surj& \overline{F}_r
\end{array}
\end{equation*}
  where $\rk_{\Ct}(\bar{\cJ}_r /\bar{\cJ}_{r+1})=\dim (\bar{J}_r 
  /\bar{J}_{r+1})={r+2\choose 2}-{r+1 \choose 2} =r+1$. 
  Suppose for instance that $r$ is odd. Then lemma~\ref{lem:f4.redi}
  implies that $P(\a)=
  \prod\limits_{-r\leq i\leq r\atop i\equiv r\pmod 2} (\a-4i)$ is the 
  characteristic polynomial of the action of $\a$ on 
  ${\bar{J}_r / \bar{J}_{r+1}}$. Therefore (and since all the
  roots are simple) the characteristic polynomial of the action of 
  $\ha$ on ${\bar{\cJ}_r / \bar{\cJ}_{r+1}}$ is $P_t(\a)=
  \prod\limits_{-r\leq i\leq r\atop i\equiv r\pmod 2} (\a-4i-f_i(t))$, for
  some $f_i(t)\in t\Ct$. This implies that $\bar{\cJ}_r / \bar{\cJ}_{r+1}= 
  \bigoplus\limits_{-r\leq i\leq r\atop i\equiv r\pmod 2} R_{r+1,i}$, 
  where $R_{r+1,i}$ is a free $\Ct$-module of rank $1$ with 
  $\ha=4i+f_i(t)$. The eigenvalue of $\hb$ on $R_{r+1,i}$ must be of the
  form $-8+O(t)$. The case $r$ even is analogous.
\end{pf}

\begin{thm}
\label{thm:f5.eHFF}
  The eigenvalues of $(\ha, \hb,\hg)$ in $\whff$ are
  $(-2t,8,0)$, $(\pm 4+2t, -8 ,0)$, $(\pm 8 \sqrt{-1}-2t, 8,0)$, $\ldots$, 
  $(\pm 4(g-1) \sqrt{-1}^g +(-1)^{g} 2t, (-1)^{g-1} 8 , 0)$.
\end{thm}

\begin{pf}
  As $\hg \cJ_{r-1} \subset \cJ_r$, one has $\hg^g \in \cJ_g$, i.e. 
  $\hg^g=0$ in $HFF_g^*$, so the only eigenvalue of $\hg$ in $HFF_g^*$, and
  hence in $\whff$, is zero. 
  To compute the eigenvalues of $\hb$ in $HFF^*_g$ we may restrict 
  to $HFF^*_g/(\hg)$, i.e.\ to every $\overline{\cF}_{g-k}$. Using 
  lemma~\ref{lem:f5.redi} recursively we find that all the eigenvalues of 
  $\hb$ in $\overline{\cF}_r$ are of the form $\pm 8+O(t)$. Thus all the eigenvalues
  of $\hb$ in $HFF^*_g$ are of the form $\pm 8+O(t)$. 

  To get the eigenvalues of $\hb$ in $\whff$, let
  us argue by contradiction. Suppose that there is an eigenvalue different
  from $\pm 8$. 
  By definition of $\whff$, there exists a vector $v=\p^w(X_1,z_1e^{tD_1}) 
  \in \whff$ such that $X=X_1\cup_Y A$ is a $4$-manifold
  with $b^+>1$, $z_1 \in \AA(X_1)$, $D_1\subset X_1$ with $\bd D_1=\SS^1$,
  $w \in H^2(X;\ZZ)$ with $w\cdot \S \equiv 1\pmod 2$, satisfying 
  $(\hb^2-64)^N v \neq 0$, for arbitrarily large $N$. Then there is a
  polynomial $P(\hb,t)=\prod (\hb +(-1)^{\e_i} 8-f_i(t))$ with
  $f_i(t)\in t\Ct$, $f_i(t)\neq 0$, $\e_i=0,1$, such that
  $$
   P(\hb, t) (\hb^2-64)^N v=0,
  $$
  for some $N\geq 0$. Substituting $v$ by $(\hb^2-64)^N v$, for some large
  $N$, we can suppose that $N=0$.
  Therefore $\Dws_X(z_1e^{tD+s\S}) \neq 0$, $\Dws_X(P(-4x,t) z_1e^{tD+s\S}) 
  =0$, with $D=D_1+\D$. As
$$
  D^w_X(z_1e^{tD+s\S}) =\frac{1}{2} 
  \left( \Dws_X(z_1e^{tD+s\S}) +\sqrt{-1}^{d_0-\deg z_1/2} 
  \Dws_X(z_1e^{\sqrt{-1}(tD+s\S)}) \right)/2,
$$
  for $d_0=d_0(X,w)=-w^2-\frac{3}{2}(1-b_1+b^+)$,
  we have $D^w_X(z_1e^{tD+s\S}) \neq 0$ and 
  $\Dws_X(Q(x,t) z_1e^{tD+s\S}) =0$, with
  $Q(x,t)=P(-4x,t)P(4x,\sqrt{-1}t)$.
  Moreover we can suppose that $z_1\in \AA(X)$ is homogeneous and that 
  all its components are orthogonal to $D$ as well as to $\S$, 
  which we express as 
  $z_1 \in\AA(<\S,D>^{\perp})$.
  Subsituting $D$ by a linear combination $aD+b\S$, $a\neq 0$, we can suppose
  that $D^2=0$, $D\in H_2(X;\ZZ)\subset H_2(X)$ and primitive, 
  $D\cdot \S \neq 0$.
  Then $D^w_X(Q(x,at) z_1e^{tD+s\S}) =0$.
  Also changing $w$ by $w+\S$ if necessary we can assume that 
  $w\cdot D\equiv 1\pmod 2$.

  At this stage, we represent $D$ by an embedded surface and invert the 
  roles of $D$ and $\S$. This corresponds to changing the metric: 
  we go from metrics giving a long neck when pulling $\S$ appart to 
  metrics giving a long neck when pulling $D$ appart. The Donaldson 
  invariants of $X$ do not change since $b^+>1$. Arguing as above,
  $D^w_X(Q'(x,s) z_1e^{tD+s\S}) =0$ for some polynomial $Q'(x,s)$. 
  This time we do not bother on whether $Q'$ is independent of $s$ 
  or not; we can take it to be just the characteristic polynomial of 
  $\hb$ acting on $HFF^*(D\x\SS^1, \SS^1)$. Now take the resultant of 
  $Q(x,at)$ and $Q'(x,s)$, which is a series $R(s,t) \neq 0$.
  Then $D^w_X(R(s,t) z_1e^{tD+s\S}) =0$ implies $D^w_X(z_1e^{tD+s\S})=0$, 
  which is a contradiction. This proves that the only eigenvalues of 
  $\hb$ are $\pm 8$.

  Finally, to compute the eigenvalues of $\ha$ we can restrict to 
  $\whff/(\hb^2-64,\seq{\hq}{1}{2g})$. 
  This is a subset of $\overline{HFF}{}^*_g=
  (HFF^*_g)_I/(\hg,\hb^2-64)$, which is computed in theorem~\ref{thm:f5.rHFF}.
  Moreover all the eigenvalues in $\overline{HFF}{}^*_g$ are indeed 
  eigenvalues of $\whff$ as all the vectors constructed in 
  proposition~\ref{prop:f5.10} come from $4$-manifolds with $b^+>1$.
  This completes the proof.
\end{pf}

\begin{rem}
  The author believes that the eigenvalues of $\whff$ are indeed all
  the eigenvalues of $HFF^*_g$. 
\end{rem}

\section{Fukaya-Floer homology $HFF_*(\Y,\d)$}
\label{sec:f6}

Now we deal with the Fukaya-Floer (co)homology
of the $3$-manifold $Y=\Y$ with the $SO(3)$-bundle
with $w_2=\PD[\SS^1] \in H^2(Y;\ZZ/2\ZZ)$ and loop
$\d \subset \S \subset \Y$ representing a primitive homology class, 
and its $\AA(\S)$-module structure. Poincar\'e duality identifies 
the Fukaya-Floer homology $HFF_*( Y,\d)$ with the Fukaya-Floer cohomology $HFF^*(Y,-\d)$.
The $\mu$ map gives an action of $\AA(\S)$ on $HFF_*( Y,\d)$. 
We shall see later that this gives in 
fact a structure of module over $HF_*(Y)$.

\subsection{The vector space $HFF_*(\Y,\d)$}
We can suppose that the basis $\{\seq{\g}{1}{2g}\}$ of
$H_1(\S;\ZZ)$ is chosen so that $[\d]=\g_1$ (recall that
$\g_i \g_{i+g}=\point$ for $1 \leq i \leq g$).
The action of $\Spz$ on $\{\g_i\}$ restricts to an action of the subgroup $\Spzr$ on
$\g_2,\ldots, \g_g,\g_{g+2},\ldots, \g_{2g}$. Any element of $\Spzr$ can
be realized by a diffeomorphism of $\Y$ fixing $\d$, hence it induces
an automorphism of $HFF_*(Y,\d)$. This gives an action of $\Spzr$ on $HFF_*(Y,\d)$.

We recall that for computing $HFF_*(Y,\d)$ there is a spectral sequence whose
$E_3$ term is $HF_*(Y) \ox \hH_*(\CP^{\infty})$, with differencital $d_3$ given by 
$$
   \mu(\d):HF_i(Y) \ox H_{2j}(\CP^{\infty}) \ar HF_{i-3}(Y) \ox
   H_{2j+2}(\CP^{\infty}).
$$
and converging to $HFF_*(Y,\d)$. The $\Spzr$ action on this $E_3$ term gives
the action on $HFF_*(Y,\d)$.
Now we can use the description of $HF_g^*=HF^*(Y)$ 
gathered in proposition~\ref{prop:f4.floer} and the fact that $\mu(\d)$ is multiplication
by $\q_1=\p^w(A,\g_1)$ to get a description of the $E_4$ term of spectral sequence.

\begin{prop}
\label{prop:f6.q1} 
  Consider $\q_1:HF^*_g \ar HF_g^*$. Then
  $\ker\q_1/\im \q_1= \bigoplus\limits_{k=0}^{g-1} \L^k_0 H^3_{\red} \ox K_{g-k}$,
  where $H^3_{\red}=<\q_2,\ldots, \q_g,\q_{g+2},\ldots,\q_{2g}>$ and $K_r={J_{r-1}/
  (J_r +\g J_{r-2})}$.
\end{prop}

\begin{pf}
The space $H^3$ has basis $\q_1, \q_2,\ldots, \q_{2g}$, so we can write 
$H^3=<\q_1,\q_{g+1}>\oplus H^3_{\red}$, where $H^3_{\red}$ is generated
by $\q_2,\ldots, \q_g,\q_{g+2},\ldots,\q_{2g}$ and it is the standard representation
of $\Spzr$ (`$\red$' stands for reduced and 
follows the notation of~\cite{genus2}). More intrinsically, we can identify
$H^3_{\red} \iso <\q_1>^{\perp}/<\q_1>$.
It is easy to check that $\L_0^k H^3$ decomposes as
$$
  \L_0^k H^3=\g'\L_0^{k-2}H^3_{\red} \oplus \left( <\q_1,\q_{g+1}>\ox \L_0^{k-1}
  H^3_{\red} \right) \oplus \L_0^k H^3_{\red}
$$
as $\Spzr$-representations, where $\g'=-g\q_1\wedge\q_{g+1}+\g$. The reader can
check this directly, noting that $\g' \in \L_0^2 H^3$, 
or otherwise see~\cite[formula (25.36)]{Fulton}.

As a shorthand, write $F_r=\Cabg/J_r=(HF_r^*)_I$.
Then proposition~\ref{prop:f4.floer} says that 
\begin{equation} \begin{array}{rcl}
  HF_g^*= \bigoplus\limits_{k=0}^{g-1} (\L_0^k H^3 \ox F_{g-k}) & = &\bigoplus\limits_{k=0}^{g-1}
  \L_0^k H^3_{\red} \ox (F_{g-k}\oplus \g'F_{g-k-2})
  \oplus  \\ & & 
  \bigoplus\limits_{k=0}^{g-1}\L_0^k H^3_{\red} \ox (<\q_1,\q_{g+1}> \ox F_{g-k-1}),
\end{array}\label{eqn:f6.dia}
\end{equation}
as $\Spzr$-representations.  
Now multiplication by $\q_1$ is $\Spzr$-equivariant and intertwines the two summands
in~\eqref{eqn:f6.dia}, i.e.
\begin{equation} 
\begin{array}{ccc}
  F_{g-k}\oplus \g'F_{g-k-2} &\stackrel{\q_1}{\ar} &<\q_1,\q_{g+1}> \ox F_{g-k-1} \\
    x \oplus \g' y &\mapsto & \q_1 \ox (x +\g y) 
\end{array}
\label{eqn:f6.uno}
\end{equation}
and
\begin{equation} 
\begin{array}{ccc}
   <\q_1,\q_{g+1}> \ox F_{g-k-1} &\stackrel{\q_1}{\ar} & F_{g-k}\oplus \g'F_{g-k-2} \\
    \q_1 \ox z &\mapsto & 0\qquad \\  
   \q_{g+1} \ox z &\mapsto & {1\over g}\g z \oplus (-{1\over g})\g' z
\end{array}
\label{eqn:f6.dos}
\end{equation} 
In~\eqref{eqn:f6.uno}, $\ker \q_1=
\{x \oplus \g'y \in F_{g-k}\oplus \g'F_{g-k-2} /  x+\g y =0 \in F_{g-k-1}\}$ 
and  $\im \q_1= \q_1 \ox F_{g-k-1}$.
In~\eqref{eqn:f6.dos}, $\ker \q_1= \q_1 \ox F_{g-k-1}$ and $\im\q_1=
\{ \g y\oplus (-\g' y) \in F_{g-k}\oplus \g'F_{g-k-2}\}$, so
$$
  \ker \q_1/\im \q_1= \bigoplus_{k=0}^{g-1} \L_0^k H^3_{\red} \ox K_{g-k},
$$
where 
\begin{eqnarray*}
K_{r} &=&
{ \{x \oplus \g'y \in F_{r}\oplus \g'F_{r-2} /  x+\g y =0 \in F_{r-1}\}   \over 
    \{ \g y\oplus (-\g' y) \in F_{r}\oplus \g'F_{r-2}\} } \iso \\
   & \iso & { \{ x \in F_{r} / x=0 \in F_{r-1}\} \over \{ \g y / y=0 \in F_{r-2}\} }= 
   {J_{r-1}/J_r  \over \g (J_{r-2}/J_r)} =
   {J_{r-1} \over J_r +\g J_{r-2}}.
\end{eqnarray*}
\end{pf}

\begin{lem}
\label{lem:f6.Kr}
  As a $\Cabg$-module, 
  $K_r=\bigoplus\limits_{-(r-1) \leq i \leq r-1 \atop i \equiv r-1 \pmod 2} R_i$, where
  $R_i$ is $1$-dimensional, $\a$ acts as $4i\sqrt{-1}$ if $i$ is even and as $4i$ if 
  $i$ is odd,
  $\b$ as $(-1)^i8$ and $\g$ as zero on $R_i$.
\end{lem}

\begin{pf}
$K_r$ is generated, as $\Cabg$-module,
by three elements $R_{r-1}^1$, $R_{r-1}^2$ and $R_{r-1}^3$, which satisfy six relations
$R_{r}^1=0$, $R_{r}^2=0$, $R_{r}^3=0$, $\g R_{r-2}^1=0$, 
$\g R_{r-2}^2=0$ and $\g R_{r-2}^3=0$. Therefore
  $$
   \left\{ \begin{array}{l} 0 = \a R_{r-1}^1 + (r-1)^2 R_{r-1}^2
   \\ 0= (\b+(-1)^{r}8) R_{r-1}^1 + {2(r-1) \over r}  R_{r-1}^3
   \\ 0 = \g  R_{r-1}^1
    \end{array} \right.
  $$
Also $R_{r-1}^3=\g R_{r-2}^1=0$. The first line allows to write $R_{r-1}^2$ in
terms of $R_{r-1}^1$, so $K_r$ is generated by an element $k_r=R_{r-1}^1$, which
satisfies $\g k_r=0$ and $(\b+(-1)^{r}8) k_r=0$. Therefore $K_r$ is a module
over $\Cabg/((\g,\b+(-1)^{r}8)+J_r)$. This is a quotient of $(HF_r^*)_I$ which has been 
computed in~\cite[proposition 20]{floer} to be 
$$
 S_r= 
  \left\{ \begin{array}{ll} \CC[\a]/((\a-16(r-1)^2)(\a-16(r-3)^2) \cdots (\a-16\cdot 1^2) ) 
  & \text{$r$ even}\\
      \CC[\a]/((\a+16(r-1)^2)(\a+16(r-3)^2) \cdots (\a+16\cdot 2^2)\a ) & \text{$r$ odd} \end{array} \right.
$$
So $K_r$ is a quotient of $S_r$, being a cyclic module over this ring. In particular
$\dim K_r \leq r$. On the other hand, if we consider the action of $\g$ in 
$F_r$,~\cite[corollary 18]{floer}
says that $\ker \g=J_{r-1}/J_r$. Moreover $\ker \g^2=J_{r-2}/J_r$, which
is proved in the same fashion. So
we can write $K_r={\ker \g \over \g \ker \g^2}$. Now 
$\dim \ker \g={r+1 \choose 2}$, $\dim \ker \g^2={r+1 \choose 2} + {r \choose 2}$. As the
action of (multiplication by) $\g$ vanishes on $\ker \g \subset \ker \g^2$, we have that
$\dim (\g\ker \g^2) \leq {r \choose 2}$. So 
$\dim (\ker \g/(\g\ker\g^2)) \geq {r+1\choose 2}-{r\choose 2}=r$, and thus $K_r$ must 
equal $S_r$.
\end{pf}

Now we are able to write down the $E_4$ term of the spectral sequence. Decompose
$\ker \q_1=\im \q_1 \oplus (\ker \q_1/\im \q_1)$, where $\im \q_1 \subset \ker \q_1$ is the
null part for the intersection pairing on $\ker \q_1$. Then 
$$
   E_4= (\im \q_1 \oplus (\ker \q_1/\im \q_1)) \x (\ker \q_1/\im \q_1) t \x 
  (\ker \q_1/\im \q_1) {t^2 \over 2!} \x \cdots
$$
So lemma~\ref{lem:f6.Kr} gives 
$$
   E_4= \im \q_1 \oplus \bigoplus_{i,k} \L^k_0 H^3_{\red} \ox R_i \ox \Ct,
$$
where $0 \leq k \leq g-1$, $-(g-k-1)\leq i \leq g-k-1$ and $i \equiv g-k-1 \pmod 2$.
We can write $E_4=\im \q_1 \oplus \tilde E_4$, where the intersection pairing vanishes on the 
first summand. In order to compute Donaldson invariants, this first summand is ineffective, so
we will ignore its behaviour through the spectral sequence, and look henceforth to 
the spectral sequence given by $\tilde E_4$.

\begin{prop}
\label{prop:f6.E4}
  The spectral sequence $\tilde E_n$, $n \geq 4$, 
  collapses at the fourth stage, i.e. $d_n=0$, for all $n \geq 4$.
\end{prop}

\begin{pf}
  There is a well-defined $\AA(\S)$-module structure in the spectral sequence, since
  it is defined at the chain level in section~\ref{sec:f3}. Also any $f \in \Spzr$
  induces $f: HFF_*(\Y, \d) \ar HFF_*(\Y,\d)$ which can be defined at the chain level
  and therefore also appears through the spectral sequence. Therefore every differential
  $d_n$ is $\Spzr$-equivariant, $\Cabg$-linear and $\Ct$-linear. 
  Now 
  $$
   \tilde E_4= \bigoplus_{i,k} \L_0^k H^3_{\red} \ox R_i \ox \Ct
  $$
  is a direct sum of inequivalent irreducible representations of $\Spzr \x \Cabg\ox \Ct$. 
  So $d_n$ has
  to send every summand to itself, and $d_n^2=0$ on it implies $d_n=0$. The proposition follows.
\end{pf}

Henceforth we will only consider
\begin{equation}
  \overline{HFF}_*(Y,\d)= 
   \bigoplus_{i,k} \L_0^k H^3_{\red} \ox R_i \ox \Ct \subset HFF_*(Y,\d),
\label{eqn:f6.lineHFF}
\end{equation}
which coincides with $HFF_*(Y,\d)/\text{null part}$. 

\subsection{The $HF^*_g$-module $HFF_*(\Y,\d)$}
In order to determine the $\AA(\S)$-module structure on $\overline{HFF}_*(Y,\d)$, we 
consider the natural cobordism between $(Y,\d) \sqcup \, (Y,\o)$ and $(Y,\d)$. It 
gives the map~\eqref{eqn:f3.3}
$$
  \cdot : HF_*(Y) \ox HFF_*(Y,\d) \ar HFF_*(Y,\d).
$$
Now for any $\p \in HFF_*(Y,\d)$, it is $\a \cdot \p=\p^w(A,2\,\S)\cdot \p=
 2 \mu(\S)({\bf 1})\cdot \p={\bf 1}\cdot 2\mu(\S)(\p) =2\mu(\S) (\p)$, with ${\bf 1}=
  \p^w(A, 1)$. Therefore
the action of $2\mu(\S)$ is multiplication by $\a$. Analogously for $\mu(\point)$ and
$\mu(\g_j)$. Therefore the $\AA(\S)$-module structure reduces to an $HF^*(Y)$-module
structure on $HFF_*(Y,\d)$, and hence on $\overline{HFF}_*(Y,\d)$.
In~\eqref{eqn:f6.lineHFF}, it is $i \equiv g-k-1 \pmod 2$, so
the action of $\mu(\g_j)$ vanishes. Therefore we have proved

\begin{thm}
\label{thm:f6.main}
  Let $Y=\Y$ and $\d\subset \S\subset Y$ a loop representing a primitive homology class. Let
  $\overline{HFF}_*(Y,\d)$ be $HFF_*(Y,\d)$ modulo its null part under the intersection
  pairing. Then $\overline{HFF}_*(Y,\d)$ is an $HF^*(Y)$-module and
\begin{equation}
  \overline{HFF}_*(Y,\d)= \bigoplus_{i,k} \L_0^k H^3_{\red} \ox R_i \ox \Ct,
\label{eqn:f6.yava2}
\end{equation}
  where $0 \leq k \leq g-1$, $-(g-k-1)\leq i \leq g-k-1$ and $i \equiv g-k-1 \pmod 2$.
  The $R_i$ are $1$-dimensional, $\a$ acts as $4i\sqrt{-1}$ if $i$ is even and as $4i$ if 
  $i$ is odd and
  $\b$ as $(-1)^i8$. The action of $\q_j$ is zero.
\end{thm}

Theorem~\ref{thm:f6.main} gives us the action of $H_*(\S)$ on $\overline{HFF}_*(Y,\d)$, but
to get a more intrinsic picture which does not need explicitly the isomorphism
$Y\iso\Y$, we need to give the action of the full $H_*(Y)$ on the Fukaya-Floer cohomology.
This is provided by the following

\begin{prop}
\label{prop:f6.act}
  Consider $\overline{HFF}_*(Y,\d)$ as given in~\eqref{eqn:f6.yava2}. Then on $R_i\ox\Ct$,
  $-4\mu(\point)$ 
  acts as $(-1)^i8$, $\mu(a)=0$ for any $a \in H_1(Y)$ and,
  for $a \in H_2(Y)$, $2\mu(a)$ is
  $4(a\cdot\SS^1)i\sqrt{-1}- 2(a \cdot \d) t$ if $i$ is even and 
  $4(a\cdot\SS^1)i+2(a \cdot \d) t$ if $i$ is odd.
\end{prop}

\begin{pf}
As $Y=\Y$ is a (trivial) circle bundle over $\S$, we may consider an automorphism of $Y$ as
a circle bundle. This is classified by an element $f \in H^1(\S;\ZZ)$, so we shall put
$\v_f : Y \ar Y$. The action in homology $\v_f: H_*(Y) \ar H_*(Y)$ is 
$\v_f(\point)=\point$, $\v_f(\g_j )=\g_j+(f[\g_j]) \SS^1$, $\v_f(\S )= \S + \PD[f]\x\SS^1$ 
and $\v_f(\a\x\SS^1 )= \a\x\SS^1$, for any $\a \in H_*(\S)$. In particular,
$$
  \d_f=\v_f(\d) =\d +n \SS^1, \qquad \text{where $n=f[\d]$.}
$$
So $\v_f:HFF_*(Y,\d) \isom HFF_*(Y,\d_f)$ and hence
\begin{equation}
  \overline{HFF}_*(Y,\d+n\SS^1)= \bigoplus_{i,k} \L_0^k H^3_{\red} \ox R_i \ox \Ct.
  \label{eqn:f6.i}
\end{equation}
Now there is a natural cobordism between $(Y,\d_f) \sqcup (Y, n\SS^1)$ and $(Y,\d_f)$, 
which gives, in the same fashion as above, an $HFF^*(Y,n\SS^1)$-module structure to 
$HFF_*(Y,\d_f)$. This goes down to a module structure over the reduced
Fukaya-Floer homology $\overline{HFF}{}^*(Y,n\SS^1)=HFF^*(Y,n\SS^1)/(\hb^2
-64,\seq{\hq}{1}{2g})$.
Corollary~\ref{cor:f5.nSS1} (and the description of
the eigenvalues of $\overline{HFF}{}^*_g$ 
given in theorem~\ref{thm:f5.rHFF}) yields
that on the summand $R_i \ox \Ct$ of $\overline{HFF}_*(Y,\d+n\SS^1)$, 
$2\mu(\S)$ must act as $4i\sqrt{-1} -2nt$ if $i$ is even and as $4i+2nt$
if $i$ is odd, $-4\mu(\point)$ as $(-1)^i8$ and $\mu(\g_j)$ as zero.
Finally we go back under the isomorphism $\v_f : Y \ar Y$. So on the summand $R_i \ox \Ct$
of~\eqref{eqn:f6.yava2}, the $\mu$-actions are as follows
$$
  \left\{
  \begin{array}{l}
   2\mu(\v_f^{-1}(\S))= 2\mu(\S - \PD[f]\x\SS^1))=
   \left\{ \begin{array}{ll}  4i\sqrt{-1} -2nt \qquad & i \text{ even}\\
              4i+2nt & i \text{ odd} \end{array} \right. \\
  -4\mu(\v_f^{-1}(\point))= -4\mu(\point)= (-1)^i 8 \vspace{2mm}\\
  \mu(\v_f^{-1}(\g_j))= \mu(\g_j -(f[\g_j]) \SS^1)=0
\end{array} \right.
$$
This implies that $\mu(\SS^1)$ acts as zero and $\mu(\g_j\x\SS^1)$ acts as 
$(-1)^{i} (\g_j\cdot \d) t$. The proposition follows.
\end{pf}

\section{Applications of Fukaya-Floer homology}
\label{sec:f7}

In this section we are going to give a number of remarkable applications from the knowledge
of the structure of the Fukaya-Floer homology groups of $\Y$. 
The author expects to extend the techniques to be able to get the
general shape of the Donaldson invariants of $4$-manifolds not of simple type with $b^+>1$.

\subsection{$4$-manifolds are of finite type}
In~\cite{Kr} it is conjectured that any $4$-manifold with $b^+>1$ is of finite type. In~\cite{froyshov}, Fr{\o}yshov gives a proof of the finite
type condition for any simply connected $4$-manifold by studying the general properties of the map 
$\mu(\point)$ on the Floer homology of $3$-manifolds.
In~\cite{wieczorek}, Wieczorek also proves the finite type condition for
simply connected $4$-manifolds by studying configurations of embedded 
spheres of negative self-intersections. Here we give a proof of the
finite type condition for arbitrary $4$-manifolds with $b^+>1$ by
using the effective Fukaya-Floer homology $\whff$.

\begin{prop}
\label{prop:f7.finite}
  Let $X$ be a $4$-manifold with $b^+> 1$ and $\S \inc X$ an embedded
  surface of self-intersection zero. Suppose there is $w \in H^2(X;\ZZ)$ with 
  $w \cdot \S \equiv 1 \pmod 2$. Then there exists $n \geq 0$ such that
  $D^w_X((x^2-4)^n z)=0$ for any $z \in \AA(X)$. 
\end{prop}

\begin{pf}
  If the genus $g$ of $\S$ is zero then the Donaldson invariants vanish 
  identically, so the statement is true with $n=0$. Suppose then that 
  $g \geq 1$. Then we split $X=X_1 \cup_Y A$, where $A$ is a small tubular 
  neighbourhood of $\S$. Let $D \in H_2(X)$ such that $D\cdot \S=1$. 
  Represent $D$ by a $2$-cycle intersecting transversely $\S$ in one positive
  point and put $D=D_1+\D$, with $D_1 \subset X_1$ and $\bd D_1 =\SS^1$. 
  Then for any $z \in \AA(X_1)$ it is $\p^w(X_1,z e^{tD_1})\in \whff$ 
  by definition~\ref{def:f5.eHFF}. By theorem~\ref{thm:f5.eHFF}
  there is some $n>0$ such that $(\hb^2-64)^n=0$ in $\whff$. Hence
  $$
   \Dws_X(z (x^2-4)^n e^{tD})=\la {1\over 16^n}(\hb^2-64)^n
  \p^w(X_1, z e^{tD_1}), \p^w(A, e^{t\D})\ra=0.
  $$
  So $D^w_X(z (x^2-4)^n D^m)=0$ for all $m\geq 0$. 
  This is equivalent to the statement.
\end{pf}

\begin{thm}
\label{thm:f7.finite}
  Let $X$ be a $4$-manifold with $b^+>1$. Then $X$ is of $w$-finite type, for
  any $w \in H^2(X;\ZZ)$.
\end{thm}

\begin{pf}
  First note that if $\tilde X=X\# \overline{\CC\PP}^2$ is the blow-up of $X$ with exceptional
  divisor $E$, then $X$ is of $w$-finite type if and only if $\tilde X$ is of $w$-finite type
  if and only if $\tilde X$ is of $(w+E)$-finite type. This is a consequence of the general
  blow-up formula~\cite{FS-bl}.
  This means that, after possibly blowing-up, we can suppose $w$ is odd. Then there exists 
  $x \in
  H_2(X;\ZZ)$ with $w\cdot x \equiv 1 \pmod 2$. As $b^+>0$, there is $y \in
  H_2(X;\ZZ)$ with $y \cdot y>0$. Consider $x'=x+2n y $ for $n$ large. 
  Then $x'\cdot x'>0$ and $w \cdot
  x'  \equiv 1 \pmod 2$. Represent $x'$ by an embedded surface $\S'$ and blow-up $X$ at 
  $N=x'\cdot x'$ points in $\S'$ to get a $4$-manifold $\tilde X=X\# N\overline{\CC\PP}^2$
  with an embedded surface $\S\subset \tilde X$ such that $\S\cdot \S=0$ and 
  $w \in H^2(X;\ZZ)\subset H^2(\tilde X;\ZZ)$ with $w\cdot\S \equiv 1 \pmod 2$. Then
  proposition~\ref{prop:f7.finite} implies that $\tilde X$ is of $w$-finite type and hence $X$
  is of $w$-finite type.
\end{pf}

\begin{prop}
\label{prop:f7.orderft}
  Let $X$ be a $4$-manifold with $b^+>1$ and containing an embedded surface $\S$ of 
  genus $g$ and self-intersection zero such that there is $w \in H^2(X;\ZZ)$ with 
  $w\cdot \S\equiv 1\pmod 2$. Then
  $X$ is of $w$-finite type of order less or equal than
  $$
   \sum_{i=1}^g \left(\left[{2g-2i \over 4}\right] +1\right),
  $$
  where $[x]$ denotes the entire part of $x$.
  If furthermore $X$ has $b_1=0$, then $X$ is of $w$-finite type of order less or equal than
\begin{equation}
   \left[{2g-2 \over 4}\right] +1.
\label{eqn:f7.orderft}
\end{equation}
\end{prop}

\begin{pf}
  The result is obvious for $g=0$. We can thus suppose $g\geq 1$.
  We only need to find the minimum $n \geq 0$ such that
  $(\hb^2-64)^n=0$ in $\whff$ (see definition~\ref{def:f5.eHFF}).
  Consider the element
  $$
  e_r=(\hb+(-1)^r 8)(\hb+(-1)^{r-1} 8)\stackrel{(r)}{\cdots} (\hb-8),
  $$
  for $1\leq r\leq g$. Using lemma~\ref{lem:f5.redi}
  we prove by induction that there are polynomials $P_r(\hb,t) \in
  \Ct[\hb]$ such that $e_rP_r =0\in \overline{\cF}_r$ and
  $P_r(\pm 8,t)\neq 0$ (indeed $P_r$ collects all the eigenvalues of 
  $\hb$ different from $\pm 8$). Thus $e_rP_r$ is a 
  multiple of $\hg$ in $\cF_r$. Now the inclusion 
  $\hg \cJ_r \subset \cJ_{r+1}$ yields that $e_rP_r \cJ_r \subset
  \cJ_{r+1}$ and, by recurrence, that $\prod_{r=1}^g e_rP_r \in \cJ_g$. 
  We conclude that $\prod_{r=1}^g e_r P_r=0$ in $HFF^*_g$. As
  $P_r$ are isomorphisms over $\whff$ by theorem~\ref{thm:f5.eHFF}, 
  we have that $\prod_{r=1}^g e_r=0$
  in $\whff$. This means that we may take $n=\sum_{i=1}^g 
  (\left[{2g-2i \over 4}\right] +1)$ to get $(\hb^2-64)^n=0$ in $\whff$.

  In the case $b_1=0$, we use that $e_gP_g$ is a multiple of $\hg$ in
  $\cF_g$. As $P_g$ is an isomorphism over $\tilde{\cF}_g$, $e_g$ is
  a multiple of $\hg$ on $\tilde{\cF}_g$. The result follows easily.
\end{pf}

\begin{rem}
\label{rem:f7.orderft}
  The bound in~\eqref{eqn:f7.orderft} is in agreement with the conjecture in~\cite{Kr}.
  Let us check some simple cases in which proposition~\ref{prop:f7.orderft} was already known 
  to hold.
  For $g=0$, we get that $X$ is of zeroth-order finite type, i.e. that the
  Donaldson invariants vanish identically. For $g=1$, we get that $X$ is of simple 
  type~\cite{thesis}~\cite{MSz}. For $g=2$ we get that $X$ is of second order finite 
  type~\cite[theorem 5.16]{thesis}. If $b_1=0$ and $g=2$, $X$ is again of simple type.
\end{rem}

\subsection{Connected sums along surfaces of $4$-manifolds with $b_1=0$}

We are going to apply the description of the Fukaya-Floer homology of section~\ref{sec:f6}
to the problem of determining the Donaldson invariants of a connected sum along a Riemann
surface of $4$-manifolds with $b_1=0$ (but not necessarily of simple type). This has been
extensively studied in~\cite{genusg}.

Let $\bar X_1$ and $\bar X_2$ be $4$-manifolds with $b_1=0$ and containing
embedded Riemann surfaces $\S=\S_i \inc \bar X_i$ of the same genus $g \geq 1$,
self-intersection zero and representing odd homology classes. 
Put $X_i$ for the complement of a small open 
tubular neighbourhood of $\S_i$ in $\bar X_i$ so that $\bar X_i=X_i \cup_Y A$, 
$X_i$ is a $4$-manifold with boundary $\bd X_i=Y=\Y$.
Let $\p:\bd X_1 \ar \overline{\bd X_2}$ be an identification 
(i.e. a bundle isomorphism)
and put $X=X(\p) = X_1 \cup_{\p} X_2 = \bar X_1 \#_{\S} \bar X_2$ 
for the connected 
sum of $\bar X_1$ and $\bar X_2$ along $\S$. 
As we are only dealing with one identification, we
may well suppose that $\p=\id$. Recall~\cite[remark 8]{genusg} 
that homology orientations of both $\bar X_i$ induce a homology orientation of $X$. 
Also choose $w_i \in H^2(X_i;\ZZ)$, $i=1,2$, and $w \in H^2(X;\ZZ)$ such 
  that $w_i \cdot \S_i \equiv 1 \pmod 2$, $w \cdot \S \equiv 1 \pmod 2$, 
  in a compatible way (i.e. the restricition of $w$ to $X_i \subset X$ 
  coincides with the restriction of $w_i$ to $X_i \subset \bar X_i$).
  Also as $b_1(X_1)=b_1(X_2)=0$ it is
  $b_1(X)=0$ and $b^+(X)>1$. Moreover there is an exact sequence~\cite[subsection 2.3.1]{thesis}
\begin{equation}
  0 \ar H_2(Y) \ar H_2(X) \ar H_2(X_1, \bd X_1) \oplus H_2(X_2, \bd X_2) \ar H_1(Y) \ar 0
\label{eqn:f7.harto}
\end{equation}

 The following result gives a strong restriction on the invariants of $X$ and 
 complements the results of~\cite{genusg}. It is also in accordance with the 
 case $g=2$ studied in~\cite{genus2}.

\begin{thm}
\label{thm:f7.harto}
  The $4$-manifold $X=\bar X_1\#_{\S} \bar X_2$ is of simple type 
  with $b_1=0$ and $b^+>1$. Let $\DD_X=e^{Q/2}\sum a_i e^{K_i}$ be its
  Donaldson series. Then for all basic classes $K_i$, we have $K_i \cdot \S 
  \equiv 2g-2 \pmod 4$.
\end{thm}

\begin{pf}
  Fix $D_S \in H_2(X)$ with $D_S|_Y=[\SS^1]\in H_1(Y)$. Now for any $\d \in H_1(\S;\ZZ)$ 
  which is primitive we consider any $D\in H_2(X)$ with $D|_Y=\d$. Represent $D+nD_S$ 
  as $D_1 +D_2$, with $D_i\subset X_i$ and $\bd D_1=\d+n\SS^1$, where $n\in\ZZ$. The Fukaya-Floer homology
  $\overline{HFF}_*(Y,\d+n\SS^1)$ has been determined in~\eqref{eqn:f6.i} and in particular
  $\hb^2-64=0$. So for any $z_1 \in \AA(X_1)$ and $z_2 \in \AA(X_2)$
  $$
   \Dws_X (z_1z_2(x^2-4) e^{t(D+nD_S)})=\la \p^w(X_1, (x^2-4)z_1 e^{tD_1}),\p^w(X_2, z_2 e^{tD_2})\ra=
  $$
  $$
   \qquad\qquad =\la {1\over 16}(\hb^2-64)\p^w(X_1, z_1 e^{tD_1}),\p^w(X_2, z_2 e^{tD_2})\ra=0.
  $$
  By continuity this implies that $\Dws_X (z (x^2-4) e^{tD})=0$ for any $D\in H_2(X)$. So
  $X$ is of $w$-simple type, and hence of simple type. 

  Now $X$ has $b_1=0$ and $b^+>1$, so we have $\DD_X=e^{Q/2}\sum a_i e^{K_i}$. Also
$$
  \p^w(X_1, e^{tD_1}) \in \overline{HFF}_*(Y,\d+n\SS^1)_I= \bigoplus_{ -(g-1)\leq i \leq g-1
  \atop i \equiv g-1 \pmod 2} R_i \ox \Ct.
$$
Put 
$$
  p(\S)=\left\{\begin{array}{ll} 
      (\S^2-(2g-2)^2)(\S^2-(2g-6)^2) \cdots (\S^2-2^2) & \text{$g$ even} \\
      (\S^2+(2g-2)^2)(\S^2-(2g-6)^2) \cdots (\S^2+4^2)\,\S  \qquad & \text{$g$ odd} 
\end{array}\right.
$$
so that  in $\overline{HFF}_*(Y,\d+n\SS^1)_I$ it is $p(\a/2+nt)=0$ if $g$ odd and 
$p(\a/2-nt)=0$ if $g$ even (see proposition~\ref{prop:f6.act}). 
Suppose for concreteness that $g$ is even (the other 
case is analogous). Then 
$$
p({\bd\over \bd s}-nt)\Dws_X (e^{t(D+nD_S)+s\S})=\Dws_X (p(\S-nt) e^{t(D+nD_S)+s\S})=
$$
$$
\qquad\qquad =\la p(\a/2-nt) \p^w(X_1,  e^{tD_1}),\p^w(X_2,e^{tD_2+s\S})\ra=0.
$$
On the other hand, as $Q(t(D+nD_S)+s\S)=Q(t(D+nD_S)) +2nts$,~\cite[proposition
12]{genus2}
implies
$$ 
  \Dws_S(e^{t(D+nD_S)+s\S})=e^{Q(t(D+nD_S))/2 +nts}\sum_{K_i\cdot\S \equiv 2\pmod 4}
  a_{i,w} e^{K_i\cdot (D+nD_S)t +(K_i\cdot \S)s} + 
$$
$$
  \qquad \qquad+e^{-Q(t(D+nD_S))/2 -nts}
  \sum_{K_i\cdot\S \equiv 0\pmod 4} a_{i,w}e^{\sqrt{-1}K_i\cdot (D+nD_S)t +\sqrt{-1}(K_i\cdot \S)s},
$$
which is a sum (over $\Ct$) 
of exponentials of the form $e^{nts+2rs}$, $-(g-1) \leq r \leq g-1$, $r \equiv 1 \pmod 2$,
and $e^{-nts+2r\sqrt{-1}s}$, $-(g-1) \leq r \leq g-1$, $r \equiv 0 \pmod 2$. 
So for $\Dws_X(e^{t(D+nD_S)+s\S})$ to be a solution of 
the ordinary differential equation $p({\bd\over \bd s}-nt)$,
the only exponentials appearing should be $e^{nts+2rs}$, with
$-(g-1) \leq r \leq g-1$, $r \equiv 1 \equiv g-1 \pmod 2$. The result follows.
\end{pf}

 From~\cite[corollary 13]{genusg}, the sum of the coefficients of all basic classes $K_i$ 
 of $X$ with $K_i\cdot \S=2r$ is zero whenever $|r|<g-1$. It is natural to expect that 
 actually these basic classes do not appear. Theorem~\ref{thm:f7.harto}
 shows that this is in fact true for $r \not{\!\! \equiv} g-1 \pmod 2$.

\subsection{Donaldson invariants of $\S_g \x \S_h$}
Our final intention is to give the Donaldson invariants of the $4$-manifold which is
given as the product of two Riemann surfaces of genus $g\geq 1$ and $h\geq 1$.
Let $S=\S_g \x \S_h$. Then $b^+=1+ 2 g h >2$, so the Donaldson invariants are well-defined.
Recall that a $4$-manifold $X$ is of $w$-strong simple type if
$D^w_X(\g z)=0$ for any $\g \in H_1(X)$, $z \in \AA(X)$, and also
$D^w_X((x^2-4) z)=0$ for any  $z \in \AA(X)$. The structure theorem of~\cite{KM}
is also valid in this case (see~\cite{otro} for a proof using Fukaya-Floer homology
groups).

\begin{prop}[\cite{KM}~\cite{otro}]
\label{prop:f7.KM}
  Let $X$ be a manifold of $w$-strong simple type for some $w$ and 
  $b^+ >1$. Then $X$ is of strong simple type and we have 
  $\DD_X^{w}= e^{Q /2} \sum (-1)^{{K_i \cdot {w} +{w}^2} \over 2} a_i \, e^{K_i}$,
  for finitely many $K_i \in H^2(X;\ZZ)$ (called basic
  classes) and rational numbers $a_i$ (the collection is empty
  when the invariants all vanish). These classes are lifts
  to integral cohomology of ${w}_2(X)$. Moreover, for any embedded
  surface
  $S \inc X$ of genus $g$, with $S^2 \geq 0$ and representing a non-torsion
  homology class, 
  one has $2g-2 \geq S^2 +|K_i \cdot S|$. 
\end{prop}

Now suppose we are in the following situation: $\bar X_1$ and $\bar X_2$ are $4$-manifolds 
containing
embedded Riemann surfaces $\S=\S_i \inc \bar X_i$ of the same genus $g \geq 1$,
self-intersection zero and representing odd elements in homology. Consider 
$X=\bar X_1 \#_{\S} \bar X_2$,
the connected sum along $\S$ (for some identification).
Suppose that $\bar X_i$ are of strong simple type and moreover that there is an
injective map
\begin{eqnarray*}
 H_2(X) &\ar& H_2(\bar X_1) \oplus H_2(\bar X_2)\\
  D &\mapsto & (\bar D_1, \bar D_2)
\end{eqnarray*}
satisfying $D^2=\bar D_1^2+\bar D_2^2$ and $D|_{X_i}=\bar D_i|_{X_i}$, $i=1,2$. Then

\begin{prop}
\label{prop:f7.nolanum}
  In the situation above $X$ is of strong simple type. Write 
  $\DD_{X_1}=e^{Q/2}\sum a_j e^{K_j}$ and\/
  $\DD_{X_2}=e^{Q/2}\sum b_k e^{L_k}$ for the Donaldson
  series for $X_1$ and $X_2$, respectively. If $g\geq 2$ then
$$
   \DD_X(e^{tD}) = e^{Q(tD)/2}(\hspace{-5mm} \sum_{K_j \cdot \S=L_k \cdot
    \S = 2g-2}\hspace{-5mm}  2^{7g-9} a_{j}b_{k} \, e^{(K_j \cdot \bar D_1 + 
    L_k \cdot \bar D_2 +2\S\cdot D)t} + 
$$
$$     + \hspace{-5mm}\sum_{K_j \cdot \S=L_k \cdot \S= -(2g-2)}
    \hspace{-5mm} (-1)^{g-1} \,2^{7g-9} a_{j}b_{k} \, e^{(K_j \cdot
    \bar D_1 +  L_k \cdot \bar D_2-2\S\cdot D)t}).
$$
  If $g=1$ then 
$$
   \DD_X(e^{tD}) = e^{Q(tD)/2} \sum_{K_j,L_k} a_{j}b_{k}
   \, e^{(K_j \cdot \bar D_1 + L_k \cdot \bar D_2)t}( \sinh (\S\cdot D)t)^2.
$$
\end{prop}

\begin{pf}
  Let us see first that $X$ is of strong simple type. Choose $w_i \in H^2(X_i;\ZZ)$, $i=1,2$, and $w \in 
  H^2(X;\ZZ)$ such that $w_i \cdot \S_i \equiv 1 \pmod 2$, $w \cdot \S \equiv 1 \pmod 2$, 
  in a compatible way. For any $D \in H_2(X)$ 
  with $D\cdot \S=1$, put $D= D_1+D_2$ with $\bar D_i=D_i+\D \subset X_i$.
  As $\bar X_1$ is of strong simple type,  
  $\Dws_{\bar X_1}((x^2-4)e^{t\bar D_1} z_s)=0$, for any $s \in \cS$, so
  $\p^w(X_1,e^{tD_1})$ is killed by $\hb^2-64$. Analogously $\p^w(X_1, e^{tD_1})$ is killed by
  $\q_i$, for $1\leq i \leq 2g$.
  Therefore $\Dws_X((x^2-4)e^{tD})=\la \p^w(X_1,(x^2-4)e^{tD_1}),\p^w(X_2,e^{tD_2})\ra=0$. 
  Analogously we see $\Dws_X(\g_i e^{tD})=0$, $1\leq i\leq 2g$. We leave to the reader the other 
  $\g \in H_1(X)$ not in the image of $H_1(\S)\ar H_1(X)$.

For the second assertion, suppose now $g\geq 2$.
Then $\p^w(X_1, e^{tD_1})$ lives in the reduced Fukaya-Floer homology 
$\overline{HFF}{}^*_g$ of subsection~\ref{subsec:rHFF}, which is found in
theorem~\ref{thm:f5.rHFF} to be isomorphic to $\CC^{2g-1}[[t]]$. Actually
it is the space $\CC^{2g-1}[[t]] \subset V[[t]]$ of~\cite[page 794]{genusg}. 
In~\cite{genusg} the intersection pairing restricted to $\CC^{2g-1}[[t]]$ is computed and then
$$
  \Dws_X(e^{tD})=\la \p^w(X_1,e^{tD_1}),\p^w(X_2,e^{tD_2})\ra
$$
is found. So the arguments in~\cite{genusg} carry over in our situation and the result
in~\cite[theorem 9]{genusg} is true for $X$. The statement follows.

The resut for $g=1$ is in~\cite[theorem 4.13]{thesis} and~\cite{MSz}.
\end{pf}

\begin{thm}
  Let $S=E \x F$  be the product of two Riemann surfaces of genus $g,h\geq 1$, i.e.
  $E=\S_g$ and $F=\S_h$. Arrange so that $h \leq g$. Then
  $S$ is of strong simple type and the Donaldson series are as follows.
  \begin{equation*} \left\{ \begin{array}{ll}
    \DD_S= 4^g e^{Q/2} \sinh^{2g-2} F \qquad & \text{if $h=1$} \\
    \DD_S= 2^{7(g-1)(h-1)+3} \sinh K & \text{if $g,h>1$, both even} \\
    \DD_S= 2^{7(g-1)(h-1)+3} \cosh K & \text{if $g,h>1$, at least one odd} 
  \end{array} \right. \end{equation*}
  where $K=K_S= (2g-2) F +(2h-2) E$ is the canonical class.
\end{thm}

\begin{pf}
  The result is a simple consequence of proposition~\ref{prop:f7.nolanum} noting 
  that $S=\S_1 \x\S_1$ 
  is of strong simple type (we leave the proof of this to the reader
  using the description of $HFF^*_1$)
  and also making use of the
  Donaldson series $\DD_S= 4 e^{Q/2}$ given in~\cite{Stipsicz}.
\end{pf}

{\em Acknowledgements.\/} 
I am grateful to the Mathematics 
Department in Universidad de M\'alaga for their hospitality and support.  
Conversations with Marcos Mari\~no, Tom Mrowka, Cliff Taubes,
Bernd Siebert, Gang Tian, Dietmar Salamon and Rogier Brussee
have been very helpful.
Special thanks to
Simon Donaldson and Ron Stern for their encouragement.


\begin{thebibliography}{10}
\bibitem{otro}
V. Mu\~noz,
Basic classes for four-manifolds not of simple type,
{\tt math.DG/9811089}.
\end{thebibliography}
\end{document}